\input ./style/arxiv-general.cfg
\documentclass[aop,MSNbibl,dvips]{arximspdf}
\makeatletter
   \@ifpackageloaded{graphicx}{}{\usepackage{graphicx}}
\makeatother

\doi{10.1214/14-AOP996}
\volume{44}
\issue{2}
\pubyear{2016}
\firstpage{1053}
\lastpage{1106}
\docsubty{FLA}

\makeatletter
\newcommand{\ds}{\displaystyle}
\newcommand{\rrvert}{\vert}
\newcommand{\rrVert}{\Vert}
\newcommand{\llvert}{\vert}
\newcommand{\llVert}{\Vert}
\newcommand{\eps}{\varepsilon}
\newcommand{\PP}{\mathbb{P}}
\newcommand{\EE}{\mathbb{E}}
\newcommand{\RR}{\mathbb{R}}
\newcommand{\TT}{\mathbb{T}}
\newcommand{\NN}{\mathbb{N}}
\renewcommand{\L}{\mathcal{L}}
\newcommand{\A}{\mathcal{A}}

\newcommand{\B}{\mathcal{B}}
\newcommand{\C}{\mathcal{C}}
\newcommand{\N}{\mathcal{N}}

\newcommand{\R}{\mathcal{R}}

\newcommand{\calE}{\mathcal{E}}
\renewcommand{\H}{\mathsf{H}}

\newcommand{\ZZ}{\mathbb{Z}}

\newcommand{\Var}{\operatorname{\mathsf{Var}}}
\newcommand{\Corr}{\operatorname{\mathsf{Corr}}}
\newcommand{\Cov}{\operatorname{\mathsf{Cov}}}
\newcommand{\II}{\mathbf{1}}
\newcommand{\G}{\mathsf{G}}
\newcommand{\V}{\mathsf{V}}
\newcommand{\E}{\mathsf{E}}
\newproclaim{defi}{Definition}[section]
\newtheorem{lemm}[defi]{Lemma}
\newtheorem{prop}[defi]{Proposition}
\newtheorem{coro}[defi]{Corollary}
\newtheorem{theo}[defi]{Theorem}
\newproclaim{rem}{Remark}
\makeatother

\begin{document}
\begin{frontmatter}

\title{Noise-stability and central limit theorems for effective
resistance of random electric networks\thanksref{T1}}
\runtitle{Stability and CLT for effective resistance}

\begin{aug}
\author{\fnms{Rapha\"el}~\snm{Rossignol}\corref{}\ead[label=e1]{raphael.rossignol@ujf-grenoble.fr}}
\runauthor{R. Rossignol}
\affiliation{Universit\'{e} Grenoble Alpes}
\address{Institut Fourier, UMR 5582\\
Universit\'{e} Grenoble Alpes\\
100 rue des maths, BP 74\\
F38402 St Martin d'H\`{e}res cedex\\
France\\
\printead{e1}}
\end{aug}
\thankstext{T1}{This work was partially supported by the ANR project GeMeCoD ANR 2011
BS01 007 01.}

%
\received{\smonth{6} \syear{2012}}
%
\revised{\smonth{6} \syear{2014}}

\begin{abstract}
We investigate the (generalized) Walsh decomposition of point-to-point
effective resistances on countable random electric networks with i.i.d.
resistances. We show that it is concentrated on low levels, and thus
point-to-point effective resistances are uniformly stable to noise. For
graphs that satisfy some homogeneity property, we show in addition that
it is concentrated on sets of small diameter. As a consequence, we
compute the right order of the variance and prove a central limit
theorem for the effective resistance through the discrete torus of side
length $n$ in $\ZZ^d$, when $n$ goes to infinity.
\end{abstract}

\end{frontmatter}

\section{Introduction}\label{sec1}

Consider a piece of conductive material whose resistivity possesses
some microscopic disorder. One way to account for this disorder is to
suppose that the material is an electric network made of tiny random
resistances. Once this model is assumed, one typically wants to
understand the behaviour of the macroscopic resistivity of the
material. To make the picture more accurate, imagine that each edge of
the lattice $\ZZ^d$ is equipped with a resistance $r(e)$ belonging to
some interval $[1,\Lambda]$ (we shall not prescribe any resistance
unit). Suppose in addition that all resistances are random, independent
and identically distributed. Our macroscopic piece of material is now
the box $\B_n=\{0,\ldots,n\}^d$, two sides of which we distinguish:
$A_n=\{x\in\B_n\mbox{ s.t. }x_1=0\}$ and $Z_n=\{x\in\B_n\mbox{ s.t.
}x_1=n\}$. The \textit{effective resistance} of the box $\B_n$ is then
defined as
\[
\R_n=\inf_{\theta}\sum_{e\in\E_n}r(e)
\theta^2(e),
\]
where the sum is over the set $\E_n$ of edges inside $\B_n$ and the
infimum is taken over all unit flows on $\E_n$ from $A_n$ to $Z_n$ (all
precise definitions are postponed until Section~\ref{sec:prel}). In the
literature, the \textit{effective conductivity} is more often the main
character. It is simply the reciprocal value of the effective
resistance and can also be defined as
%
\begin{equation}
\label{eq:conductancemin}
\C_n=\inf_v\sum
_{e\in\E_n}c(e) \bigl(dv(e)\bigr)^2,
\end{equation}
where the infimum is over all functions $v$ on $\B_n$ having value $0$
on $A_n$ and $1$ on $Z_n$, $c(e)=1/r(e)$ is the conductance of edge
$e$, and $dv(e):=v(e_-)-v(e_+)$ is the difference of $v$ along edge
$e$. The unique minimizer in the definition of $\C_n$ is the function
that is $0$ on $A_n$ and $1$ on $Z_n$ and is discrete harmonic on
$\B_n\setminus\{A_n,Z_n\}$. It is worth mentioning that the setting above
is also relevant to describe the pressure field of a fluid through a
weakly porous medium when the circulation of the fluid can be modelled
with Poiseuille's law. The central problem is now to understand the
asymptotic behaviour of $\R_n$ (or, equivalently, of $\C_n$) as $n$
goes to infinity.

A first step in this direction was accomplished in the setting of
stochastic homogenization
theory (cf. \cite{ZhikovKozlovOleinik94},
Chapter~7). It is shown in \cite{Kozlov86}, Section~3, that a law of
large numbers holds (see also \cite{PapanicolaouVaradhan81,Kunnemann83} and \cite{BoivinDepauw03} for related results). Namely,
there is some positive constant $\mu$ such that
\[
\frac{1}{n^{d-2}}\C_n\mathop{\longrightarrow}_{n\to\infty}^{\mathrm{a.s.}}
\mu \quad\mbox{that is}\quad n^{d-2}\R_n\mathop{
\longrightarrow}_{n\to\infty}^{\mathrm{a.s.}}\frac{1}{\mu}.
\]
To understand\vspace*{1pt} the scaling, notice that the function $v_{\mathrm{hom}}\dvtx x\mapsto
x_1/n$ gives an upper bound of order $n^{d-2}$ on the value of $\C_n$,
and a flow $\theta$ satisfying\vspace*{1pt} $\theta(e)=\frac
{1}{(n+1)^{d-1}}dv_{\mathrm{hom}}$ gives an upper bound of order $n^{2-d}$ on
the value of $\R_n$.


A second step is to understand the fluctuations of $\C_n$ and $\R_n$.
If the optimal function in the definition of $\C_n$ were $v_{\mathrm{hom}}$,
then $\C_n$ would merely be a sum of $\Theta(n^d)$ i.i.d. random
variables, each of variance $\Theta(n^{-4})$. The variance of $\C_n$
would thus be of order $\Theta(n^{d-4})$, and that of $\R_n$ of order
$\Theta(n^{4-3d})$. A lower bound of this order was given by Wehr (cf.
\cite{Wehr97}) under some technical assumptions (see also Section~\ref{sec:lowlevels} below). More recently, an upper bound of the same order
was obtained by \cite{GloriaOttoVariance} for a different, but closely
related quantity. We shall present in more details the work of Gloria
and Otto at the end of this introduction.

The main purpose of the present paper is to derive the right order of
the variances of $\C_n$ and $\R_n$ and in addition to make a step
further in the understanding of their fluctuations by deriving Gaussian
central limit theorems for these quantities. However, for technical
reasons we shall only be able to do this in a translation invariant
setting, namely for the effective resistance through the torus; cf.
Theorem~\ref{theo:CLTTorus}. This is the main result of the article.
Our approach to obtain this result is however quite general and not
restricted to graphs like $\ZZ^d$. Indeed, we shall study in
Section~\ref{sec:ES} \textit{the generalized Walsh decomposition} of
point-to-point effective resistance on general infinite networks. This
decomposition, sometimes called the Efron--Stein decomposition, is an
extension of the Fourier--Walsh decomposition on the discrete cube and
is related to a notion of noise sensitivity introduced in~\cite{BenjaminiKalaiSchramm99}. From now on, we shall drop the term
``generalized'' for simplicity. The Walsh decomposition of a square
integrable function $f$ of the resistances reads
\[
f=\sum_Sf_S,
\]
where the sum runs over all finite subsets of the set of edges, $f_S$
is a function of $(r(e))_{e\in S}$ for any $S$ and $f_S$ is orthogonal
to $f_{S'}$ whenever $S\neq S'$. This decomposition has two immediate
interesting features. First, the variance of $f$ may be expressed as
\[
\Var(f)=\sum_{S\neq\varnothing}\|f_S
\|_2^2.
\]
Second, if $f^\eps$ is obtained from $f$ by resampling independently
each input $r(e)$ with probability $\eps$, the correlation between $f$
and its $\eps$-noised version $f^{\eps}$ equals (see Proposition~\ref{prop:noise} for further details)
\[
\Corr\bigl(f,f^{\eps}\bigr)=\frac{\sum_{S\neq\varnothing}(1-\eps)^{|S|}\|f_S\|
^2}{\sum_{S\neq\varnothing}\|f_S\|^2}.
\]
Thus, if $f$ is nonconstant one may associate a ``spectral probability
measure'' $Q_f$ to $f$ on the set of nonempty finite subsets of the
set of edges:
\[
Q_f(S)=\frac{\|f_S\|_2^2}{\Var(f)},
\]
and we see that understanding the distribution of $|S|$ under $Q_f$
allows to control the noise-sensitivity of $f$. Our first result,
Theorem~\ref{theo:EStight}, shows that the second moment of $|S|$ under
$Q_f$ when $f$ is a point-to-point effective resistance, is bounded
above by a constant depending only on $\Lambda$. This implies, loosely
speaking, that the Walsh decomposition of the effective resistance is
always concentrated on low levels. More precisely, consider the set of
distributions of $|S|$ under $Q_f$ when $f$ runs over all possible
point-to-point resistances on graphs equipped with independent
resistances in $[1,\Lambda]$ with $\Lambda$ fixed. Then our first
result implies that this set of probability measures is tight. It
implies also that effective resistances are always uniformly stable to
noise in the sense of \cite{BenjaminiKalaiSchramm99} (cf.
Corollary~\ref{coro:noise}) and that the Efron--Stein inequality is
always sharp for
estimating the variance of the effective resistance; cf.
Corollary~\ref{coro:EStight}. Then we shall improve this result on a class of graphs
which possess some homogeneity property. These graphs that we shall
qualify as having \textit{homogeneous currents} contain all
quasi-transitive graphs; cf. Corollary~\ref{coro:transitivehomogeneous}. On those graphs we shall show that the
Walsh decomposition is, loosely speaking, concentrated on sets of small
diameter. This is the key to obtain a central limit theorem since the
sets of resistances with bounded diameter exhibit only finite range
dependence; cf. Section~\ref{sec:TLC}.

In Section~\ref{sec:tori}, we shall adapt this general approach to the
effective resistance through the discrete torus, deriving the optimal
variance estimate and the Gaussian central limit theorem already mentioned.

We end this introduction by giving more details on the work \cite
{GloriaOttoVariance}, and comparing our results to theirs. Their work
is close to the homogenization theory framework. Consider the discrete
elliptic differential operator $d^*(cd(\cdot))$ corresponding to random,
translation invariant and ergodic conductances $c=(c(e))_{e\in\EE^d}$
on $\ZZ^d$. Precise definitions of $d^*$ and $d$ are given in
Section~\ref{sec:prel_current}, but let us just mention that it gives,
for a function $v$ on $\ZZ^d$:
\[
d^*(cdv) (x)=\sum_{y\sim x}c(x,y) \bigl(v(x)-v(y)\bigr).
\]
Then, using the words of \cite{GloriaOttoVariance}, homogenization
theory (namely \cite{Kunnemann83}) shows that there exists a constant
matrix $A$ such that the solution operator of $\nabla^*(A\nabla(\cdot))$
describes the large scale behaviour of the solutions operator of
$d^*(cd(\cdot))$. Furthermore, $A$ can be characterized by the so called
corrector: for any $\xi$ in $\RR^d$, there exists a unique function
$\phi_\xi$ on $\ZZ^d$ (which is a function also of the conductances)
such that $\nabla\phi_\xi$ is stationary, $\phi_\xi(0)=0$, $\EE(\nabla
\phi_\xi)=0$ and\vspace*{1pt} such that $g_\xi\dvtx x\mapsto\xi.x+\phi_\xi(x)$ is
discrete harmonic for $d^*(cd(\cdot))$ on $\ZZ^d$. Then $A$ is
characterized by
\[
\xi.A\xi=\EE [\nabla g_{\xi}.A\nabla g_{\xi}].
\]
When the conductances are i.i.d., $A$ equals $\mu$ times the identity
matrix, and the constant $\mu$ is the same as in the law of large
numbers of $\C_n$ stated above. We shall fix $\xi=(1,0,\ldots,0)$ in
the sequel. When one is interested in computing $\mu$, Gloria and Otto
remark that the preceding characterization is not computationally
tractable. Thus, one has to find a way to efficiently \textit{estimate}
$\mu$. The quantity $\C_n/n^{d-2}$ is therefore a reasonable estimator
for $\mu$.
Putting aside for a moment the problem of controlling its
bias, this is where the knowledge of its variance, and even of a
central limit theorem, may be useful. Unfortunately, $\C_n$ lacks
stationarity. This is a handicap for error analysis, as Gloria and Otto
noticed in \cite{GloriaOttoVariance} for a quantity very similar to
$\C_n$. Next, they introduce a stationary approximation of the voltage,
namely $\phi_T$ solving
\[
\frac{1}{T}\phi_T+d^*\bigl(c(\xi+d\phi_T)
\bigr)=0\qquad \mbox{in } \ZZ^d.
\]
Let $\eta_L$ be an averaging cutoff function with support in $(0,n)^d$
(and some extra regularity condition). When $T$ is large with respect
to $n$, they show that the quantity
\[
\A_n:=\sum_{e\in\EE^d}c(e) (\xi+d
\phi_T) (e)^2\eta_L(e),
\]
is a good proxy for $n^{2-d}\C_n$ and, furthermore, they show that the
variance of $\A_n$ is of order at most $n^{-d}$, with some extra
polylogarithmic factor in $T$ for $d=2$. This order coincides with the
variance order conjectured above for $n^{2-d}\C_n$.

What we shall obtain in Theorem~\ref{theo:CLTTorus} is an optimal
variance estimate and a central limit theorem for the effective
conductance on the discrete torus of length $n$ when $n$ goes to
infinity. To compare our results to those obtained by Gloria and Otto,
we shall say that the precise quantity that we analyse is practically
computable and stationary. Furthermore, in some sense, the discrete
tori converge to~$\ZZ^d$ better than the discrete cubes since they
avoid boundary effects. Thus, the effective conductance on the torus
may be a better estimator of $\mu$. Notice, however, that the
convergence of the normalized effective conductance to $\mu$ is known
(cf. \cite{CaputoIoffe03} and \cite{Owhadi03}) but not the rate of
convergence. It would be interesting to investigate this rate, for
instance, in the spirit of \cite{GloriaOttoError}. Second, our method
works the same way whether $d=2$ or not, and this is an advantage over
Gloria and Otto's result, which makes a distinction between the two.
Finally, the fact that we obtain a central limit theorem is really a
step forward compared to \cite{GloriaOttoVariance} which only obtains
variance estimates. On the other hand, Gloria and Otto obtain other
interesting results, that we do not get by our method, notably
concerning the integrability of the corrector itself (Proposition~2.1
in \cite{GloriaOttoVariance}).

After this paper was submitted for publication, we learned the
existence of two preprints which address essentially the same question.
Nolen \cite{Nolen11} defines a continuous version of the effective
conductance on the torus (but with discrete randomness) and shows a
Gaussian approximation. He uses essentially two arguments: a~second-order Poincar\'{e} inequality due to Chatterjee \cite{ChatterjeeSecondOrder},
and the results of Gloria and Otto on the
boundedness in $L^p$ of the corrector. The drawbacks of this approach
are twofold. First, the bound obtained by Nolen in dimension 2 is
suboptimal, because in dimension 2, integrability results of Gloria and
Otto are weaker. Then the use of Chatterjee's inequality forces the
elliptic conductances to have a special form of distribution (notably,
it must be absolutely continuous with respect to the Lebesgue measure).
In return, Nolen obtains a bound on the variation distance between the
normalized effective conductance and the standard Gaussian
distribution, which is of course a stronger conclusion than ours. The
other preprint is by Biskup, Salvi and Wolff \cite{BiskupSalviWolff12}.
It shows a central limit theorem for the effective conductance on the
grid with linear boundary condition. One serious limitation of their
approach is that they require a small ellipticity contrast (i.e.,
$\Lambda$ close enough to 1 in our setting). On the other hand, this
paper has the advantage of giving an asymptotic equivalent of the
variance of the effective conductance. 

\section{Preliminaries}
\label{sec:prel}
\subsection{Effective resistance and minimal current}
\label{sec:prel_current}
An excellent reference for background on electric networks is the book
\cite{LyonsPeres2011}, Chapters~2 and 9 and we shall try to stick to
its notation.

In the sequel, $\G=(\V,\E)$ will be a \textit{countable, locally finite,
oriented, symmetric and connected graph}. Symmetric means that $\E$ is
a symmetric subset of $\V^2$, that is, each edge of $\G$ occurs with
both orientations in $\E$, countable means here that both $\V$ and $\E$
are at most countable and locally finite means that every vertex has
finite degree. When $e\in\E$, we let $e_-$ denote the tail of $e$ and
$e_+$ its head, we denote by $-e:=(e_+,e_-)$ the edge $e$ with reversed
direction and let $\E_{1/2}$ be a subset of $\E$ such that for every
edge $e$, exactly one of $e$ and $-e$ belongs to $\E_{1/2}$.

For every collection $r\in(0,\infty)^{\E_{1/2}}$, one may define the
\textit{electric network} $(\G,r)$: it must be understood as a resistive
network, where each edge $e$ is a resistor with resistance $r(e)$. We
shall sometimes use the notation $c(e)$ to denote the \textit{conductance} of edge $e$, that is, $c(e)=1/r(e)$. We define \textit{the
co-boundary operator} $d$ from $\RR^\V$ to $\RR^\E$ by
\[
dv(e)=v(e_-)-v(e_+),
\]
and \textit{the boundary operator} $d^*$ from $\RR^\E$ to $\RR^\V$ by
\[
d^*\theta(x)=\sum_{e_-=x}\theta(e).
\]
Notice that $dv$ plays the role of a gradient and $d^*\theta$ the role
of a divergence.

For a fixed collection $r$, we define $\ell^2_-(\E,r)$ as the Hilbert
space of \textit{antisymmetric functions on the edges that have bounded energy}:
\[
\ell^2_-(\E,r)= \bigl\{\theta\in\RR^\E\mbox{ s.t. }
\calE_r(\theta)<\infty \mbox{ and }\forall e\in\E, \theta(e)=-
\theta(-e) \bigr\},
\]
where
\[
\calE_r(\theta):=\sum_{e\in\E_{1/2}}r(e)
\theta^2(e),
\]
endowed with the scalar product:
\[
\bigl(\theta,\theta'\bigr)_r=\sum
_{e\in\E_{1/2}}r(e)\theta(e)\theta'(e).
\]
We shall denote by $\|\theta\|_r:=(\theta,\theta)_r^{1/2}$ the norm
associated to this scalar product. Thus, $\calE_r(\theta)$ is the
square of the norm in $\ell^2_-(\E,r)$ of $\theta$ and it is called the
\textit{energy} of $\theta$. In the main part of the present paper (from
Section~\ref{subsec:homogeneous} on), we shall be interested in \textit{elliptic networks}, that is, networks $(\G,r)$ for which there is a
finite constant $\Lambda\geq1$ such that $r\in[1,\Lambda]^\E$. \textit{In the whole article}, $C(\Lambda)$
[\textit{resp.}, $C(\Lambda,\G)$] \textit{will denote
a constant}, \textit{depending only on} $\Lambda$ (\textit{resp.}, \textit{on} $\Lambda$ \textit{and} $\G$),
\textit{that may vary from time to time}.  Of course, all the sets $\ell^2_-(\E
,r)$ for $r\in[1,\Lambda]^{\E_{1/2}}$ are the same, and we shall
define this common set as $\ell^2_-(\E)$:
%
\begin{equation}
\label{eq:defellE}
\ell^2_-(\E):= \biggl\{\theta\in\RR^\E
\mbox{ s.t. }\sum_{e\in\E
_{1/2}}\theta^2(e)<\infty
\mbox{ and }\forall e\in\E, \theta(e)=-\theta (-e) \biggr\}.
\end{equation}

Let $I$ be a nonnegative real number, and $u$ and $v$ two distinct
vertices of $\G$. A~member $\theta$ of $\ell^2_-(\E,r)$ is called \textit{a flow of intensity $I$ from $u$ to $v$} if:
%
\begin{equation}
\label{eq:flow}
d^*\theta=I(\II_u-\II_v).
\end{equation}
This means $\theta$ satisfies the node law on the network, except at
$u$ where a net flow of value $I$ enters the network, and at $v$ where
a net flow of value $I$ leaves the network. When $\theta$ is a flow
from $u$ to $v$, we say that $\theta$ is a unit flow if its intensity
is 1. Among flows, some are particular important: the \textit{currents},
which satisfy Kirchhoff's cycle law as stated precisely in the
following definition.

\begin{defi}
\label{defi:Kirchhoff}
For any $e\in\E$, let $\chi_e=\II_{\{e\}}-\II_{\{-e\}}$ denote the
unit flow along $e$.
A \textit{current} $i\in\ell_-^2(\E,r)$ from $u$ to
$v$ is a flow from $u$ to $v$ which satisfies \textit{Kirchhoff's cycle
law}: if $e_1,\ldots,e_n$ is an oriented cycle in $\G$, then
\[
\Biggl(\sum_{i=1}^n\chi^{e_k},i
\Biggr)_r=0.
\]
\end{defi}

Currents are the flows which derive from a potential: if $i$ is a
current, there exists a function $v$ on $\V$ such that $r(e)i(e)=dv(e)$
for any edge $e\in\E$.

We may now define the \textit{effective resistance between two points $u$
and $v$} on the network $(\G,r)$ as the minimal energy of a unit flow
between $u$ and $v$:
\[
\R^{u,v}(r):=\inf \biggl\{\sum_{e\in\E_{1/2}}r(e)
\theta^2(e)\mbox{ s.t. }\theta\in\ell_-(\E,r)\mbox{ is a unit flow
from }u \mbox{ to }v \biggr\}.
\]
Since one minimizes a Hilbert norm on a nonempty closed convex set
(non\-empty because it contains the flows induced by the paths from $u$
to $v$), the infimum above is attained by a unique flow (cf.
Proposition~9.2 of \cite{LyonsPeres2011}). It turns out that this flow
has the additional property of being a current. It is called \textit{the
minimal unit current from $u$ to $v$} and we shall denote it by
$i^{u,v}_r$. The term \textit{minimal} stems from the fact that it
minimizes the energy among all unit currents from $u$ to $v$.

It is important to notice that currents from $u$ to $v$ of prescribed
intensity may or may not be unique depending on the particular network;
see Chapter~9 in \cite{LyonsPeres2011}. On finite networks, however, it
is well known that currents are unique (see, e.g., Chapter~2 in
\cite{LyonsPeres2011}). A useful fact about minimal currents is that
they are limits of currents on finite graphs. Let us be more precise.
Let $(\G_n)_{n\geq0}$ be a sequence of finite subgraphs of $\G$ that
\textit{exhausts} $\G$, that is, such that $\G_n\subset\G_{n+1}$ and such
that $\G=\bigcup_{n\geq0} \G_n$. Suppose that $u$ and $v$ belong to $\G
_0$ and denote by $\G_n^W$ the ``wired'' network obtained from $\G$ by
identifying all vertices outside $\G_n$ as a single vertex. Notice that
one may identify the edges of $\G_n$ and $\G_n^W$ as subsets of $\E$.
Let $i_{r,n}^W$ be the (unique) current from $u$ to $v$ on $\G_n^W$ and
see it as an element of $\ell_-(\E,r)$ by putting zero flow on edges
not in $\G_n^W$. Then $i_{r,n}^W$ converges in $\ell_-(\E,r)$ (and thus
pointwise) as $n$ goes to infinity (see Proposition~9.2 in \cite
{LyonsPeres2011}).

We finish this section with a useful lemma: the absolute value of a
minimal unit current is at most one on any edge. This is intuitively
clear since it must carry a unit mass from $u$ to $v$ and also minimize
the energy.

\begin{lemm}
\label{lem:boundcurrent}
For any distinct vertices $u$ and $v$ on a network $(\G,r)$, and any
edge $e$,
\[
\bigl|i_{r}^{u,v}(e)\bigr|\leq1.
\]
\end{lemm}

\begin{pf}
Suppose first that $\G$ is finite. Let $e$ be any edge of $G$ and
suppose, without loss of generality, that $i_r(e)> 0$. Let $f$ denote
the voltage associated to $i_r$ with value zero at $v$ (see Chapter~2
in \cite{LyonsPeres2011}). It satisfies, for any edge $e'$:
\[
df\bigl(e'\bigr)=i_r\bigl(e'\bigr),
\]
notably $f$ is discrete harmonic on $\V\setminus\{u,v\}$, that is,
\[
\forall x\in\V\setminus\{u,v\}, \qquad \sum_{y\sim x}c(x,y)
\bigl(f(x)-f(y)\bigr)=0,
\]
and $f(u)>0$ (see, e.g., equation~(2.3) in \cite{LyonsPeres2011}).
Furthermore, it satisfies the maximum principle on
$\V\setminus\{u,v\}$ (see Section~2.1 in \cite{LyonsPeres2011}): for
any $W\subset\V\setminus\{u,v\}$ let $\partial W$ be the set of
vertices which are adjacent to a vertew in $W$. Then the maximum and
the minimum of $f$ on $\overline{W}=W\cup \partial W$ are attained on
$\partial W$. Now, consider the set
\[
A=\bigl\{x\in\G\mbox{ s.t. }f(x)> f(e_+)\bigr\}.
\]
$A$ is a connected set of vertices containing $u$ and $e_-$, and not
containing $v$ nor $e_+$. Indeed, $A$ clearly contains $u$ and $e_-$.
Furthermore, if $A$ had a connected component $W$ not containing $u$,
then from the maximum principle, the maximum of $f$ on $\overline{W}$
would be obtained at some $x\in\partial W$, showing that there is
$y\in W$ such that $x\sim y$ and $f(x)\geq f(y)$, but then $x$ would be
in $A$, contradicting the fact that $W$ is a connected component of $A$.

Let $\Pi$ be the set of edges with the tail in $A$ and the head in
$A^c$. Thus, $e$ belongs to $\Pi$. Because of the node law,
\[
\sum_{x\in A}d^*i_r(x)=d^*i_r(u)=1,
\]
and on the other hand,
\begin{eqnarray*}
\sum_{x\in A}d^*i_r(x)&=&\sum
_{x\in A}\sum_{e'\in\E}i_r
\bigl(e'\bigr)\II _{e'_-=x}
\\
&=&\sum_{e'\in\E}i_r\bigl(e'
\bigr)\sum_{x\in A}\II_{e'_-=x}
\\
&=&\sum_{e'\in\E}i_r\bigl(e'
\bigr)\II_{e'_-\in A}
\\
&=&\sum_{e'\in\Pi}i_r\bigl(e'
\bigr),
\end{eqnarray*}
since $i_r$ is antisymmetric, but of course, $i_r(e')\geq0$ for any
$e'$ in $A$. Thus,
\[
0<i_r(e)\leq\sum_{e'\in\Pi}i_r
\bigl(e'\bigr)=1.
\]
This shows the result on finite graphs. It implies the general result
since the minimal unit current between $u$ and $v$ is the pointwise
limit of a sequence of minimal unit currents between $u$ and $v$ on
finite graphs.
\end{pf}

\subsection{Partial derivatives of the effective resistance and the
minimal current}
The functions $r\mapsto i^{u,v}_r$ and $r\mapsto\R^{u,v}(r)$ are
smooth functions, as the next lemma shows. In the sequel, $\partial_ef$
denotes the partial derivative with respect to $r(e)$ of a function $f$
on $(0,\infty)^{\E_{1/2}}$ and $\partial^2_{e,e'}f$ denotes $\partial
_e\,\partial_{e'}f$.

\begin{lemm}
\label{lem:regularite1}
The functions $r\mapsto i^{u,v}_r(e)$, for any edge $e$, and $r\mapsto
\R^{u,v}(r)$ admit partial derivatives of all orders. In addition, for
any distinct vertices $u,v$ and edges $e$, $e'$:
\begin{longlist}[(iii)]
\item[(i)] $\forall e'\neq e, \partial_{e'}
i^{u,v}_r(e)=\frac{i^{u,v}_r(e')}{r(e')}i^{e'}_r(e)=\frac
{i^{u,v}_r(e')}{r(e)}i^{e}_r(e') $.
\item[(ii)]
$\forall e, \partial_e i^{u,v}_r(e)=\frac
{i^{u,v}_r(e)}{r(e)}(i^{e}_r(e)-1) $.
\item[(iii)]
$\forall e, \partial_e \R
^{u,v}(r)=(i^{u,v}_r(e))^2 $.
\end{longlist}
\end{lemm}

\begin{pf}
Let us first suppose that $\G$ is finite. Then it is well known that
$i^{u,v}_r(e)$ and $\R^{u,v}(r)$ are rational functions of $r$ with no
positive pole. See, for instance, \cite{Bollobas98}, Theorem~2, page~46.
The idea goes back to Kirchhoff (see \cite{Kirchhoff58} for an english
translation of the original paper). The fact that $\partial_e \R
^{u,v}(r)=(i^{u,v}_r(e))^2$ is also known; cf., for instance, \cite
{LyonsPeres2011}, Exercise 2.69. One easy way to see it is as follows.
Let $r'$ be a collection of resistance differing from $r$ only on edge
$e$. Then, using the minimality of $i^{u,v}_r$,
\[
\calE_{r'}\bigl(i^{u,v}_{r'}\bigr)-
\calE_{r}\bigl(i^{u,v}_{r'}\bigr)\leq\bigl(
\R^{u,v}\bigl(r'\bigr)-\R ^{u,v}(r)\bigr)\leq
\calE_{r'}\bigl(i^{u,v}_{r}\bigr)-
\calE_{r}\bigl(i^{u,v}_{r}\bigr),
\]
and thus
%
\begin{equation}
\label{eq:majoDeltaeR2sided} \hspace*{3pt}\quad\bigl(r'(e)-r(e)\bigr) \bigl(i_{r'}^{u,v}(e)
\bigr)^2\leq\bigl(\R^{u,v}\bigl(r'\bigr)-
\R^{u,v}(r)\bigr)\leq \bigl(r'(e)-r(e)\bigr)
\bigl(i_{r}^{u,v}(e)\bigr)^2.
\end{equation}
Letting $r'(e)$ go to $r(e)$ shows that $\partial_e \R
^{u,v}(r)=(i^{u,v}_r(e))^2$.

To compute the partial derivatives of $r\mapsto i^{u,v}_r(e)$, let us
differentiate the flow condition~(\ref{eq:flow}) and Kirchhoff's cycle
law of Definition~\ref{defi:Kirchhoff} with respect to $r(e')$. We obtain
\[
\forall x\in\V,\qquad  d^* \bigl[\partial_{e'}i^{u,v}_r
\bigr](x)=0,
\]
and for every cycle $\gamma$ on $\G$,
\[
\sum_{e\in\gamma} \bigl[r(e)\,\partial_{e'}i^{u,v}_r(e)+i^{u,v}_r
\bigl(e'\bigr)\chi _{e'}(e) \bigr]=0.
\]
Thus, if one defines
\[
j(e)=\partial_{e'}i^{u,v}_r(e)+
\frac{i^{u,v}_r(e')}{r(e')}\chi_{e'}(e),
\]
we get
\[
\forall x\notin\bigl\{e_-',e_+'\bigr\},\qquad  d^*j(x)=0,
\]
and, for every cycle $\gamma$ on $\G$,
\[
\sum_{e\in\gamma} r(e)j(e)=0.
\]
Thus, $j$ is a current from $e_-$ to $e_+$ or from $e_+$ to $e_-$,
depending on the sign of $d^*j(e'_-)$. Its intensity is deduced from
\[
d^*j\bigl(e'_-\bigr)=d^* \bigl[\partial_{e'}i^{u,v}_r
\bigr](e_-)+\frac
{i^{u,v}_r(e')}{r(e')}=\frac{i^{u,v}_r(e')}{r(e')}.
\]
Thus, from the unicity of currents on finite graphs, one gets
\[
j=\frac{i^{u,v}_r(e')}{r(e')}i^{e'}_r.
\]
Consequently,
\[
\forall e'\neq e, \qquad \partial_{e'} i^{u,v}_r(e)=
\frac
{i^{u,v}_r(e')}{r(e')}i^{e'}_r(e)
\]
and
\[
\forall e',\qquad \partial_{e'} i^{u,v}_r
\bigl(e'\bigr)=\frac
{i^{u,v}_r(e')}{r(e')}j\bigl(e'\bigr)=
\frac
{i^{u,v}_r(e')}{r(e')}\bigl(i^{e'}_r\bigl(e'
\bigr)-1\bigr).
\]
Finally, for any $e\neq e'$,
\[
\partial^2_{e',e} \R^{u,v}(r)=\partial
_{e'}\bigl(i^{u,v}_r(e)\bigr)^2=2i^{u,v}_r(e)
\frac
{i^{u,v}_r(e')}{r(e')}i^{e'}_r(e).
\]
But since $\partial^2_{e',e} \R^{u,v}(r)=\partial^2_{e,e'} \R
^{u,v}(r)$, we obtain
\[
i^{u,v}_r(e)\frac{i^{u,v}_r(e')}{r(e')}i^{e'}_r(e)=i^{u,v}_r
\bigl(e'\bigr)\frac
{i^{u,v}_r(e)}{r(e)}i^{e}_r
\bigl(e'\bigr).
\]
Now, take $(u,v)=e$ and $(u,v)=e'$ and notice that $i^{e}_r(e)$ and
$i^{e'}_r(e')$ are always different from zero. We obtain
\[
\frac{i_r^e(e')i_r^{e'}(e)}{r(e')}=\frac{(i_r^e(e'))^2}{r(e)}
\]
and
\[
\frac{i_r^e(e')i_r^{e'}(e)}{r(e)}=\frac{(i_r^{e'}(e))^2}{r(e')}.
\]
Thus, one deduces that $i_r^e(e')=0$ if and only if $i_r^{e'}(e)=0$ and
in any case,
%
\begin{equation}
\label{eq:reciprocity}
\frac{i_r^e(e')}{r(e)}=\frac{i_r^{e'}(e)}{r(e')}.
\end{equation}
This last relation is called the \textit{reciprocity law}.
See \cite{LyonsPeres2011}, Chapter~2 for another proof. This concludes the proof
of the lemma on finite graphs.

Now, let $\G$ be infinite, $r$ belong to $(0,\infty)^{\E_{1/2}}$ and
$e$, $e'$ in $\E_{1/2}$.
As explained in Section~\ref{sec:prel_current}, for any $u$ and $v$, $i_r^{u,v}$ is the limit, in
$\ell^2_-(\E,r)$ of a sequence $i_{r,n}^W$ of unit currents from $u$ to
$v$ on ``wired'' finite graphs $G_n^W$. Notably, $i_r^{u,v}(e)$ is the
pointwise limit of $i_{r,n}^W(e)$. From the formulas of the derivatives
on finite graphs and Lemma~\ref{lem:boundcurrent}, one sees that
$r(e')\mapsto i_{r,n}^W(e)$ form\vspace*{1pt} an equi-continuous family of functions
on any compact interval $I$ of $(0,\infty)$. It follows from
Arzela--Ascoli's theorem that the convergence of $i_{r,n}^{u,v}(e)$ to
$i_r^{u,v}(e)$ is uniform when $r(e')$ runs over~$I$. Notably, this
implies the continuity of $r(e')\mapsto i_{r}^{u,v}(e)$ for any $e$,
$e'$. Then, from the formulas (i) and (ii) of the derivative
$\partial_{e'}i_{r,n}^{u,v}(e)$ one sees that the derivative itself
converges uniformly when $r(e')$ runs over $I$. Then $r(e')\mapsto
i_r^{u,v}(e)$ is differentiable on $(0,\infty)$ and its derivative is
the limit of the derivatives $\partial_{e'}i_{r,n}^{u,v}(e)$. This
shows the formulas for $\partial_{e'}i_{r,n}^{u,v}(e)$. Formula (iii)
is then a consequence of (\ref{eq:majoDeltaeR2sided}), since
$r(e)\mapsto i_{r}^{u,v}(e)$ is continuous.
\end{pf}

\begin{rem}
A similar formula, relating the partial derivative of the voltage drop
through $e$ with respect to the conductance $c(e')$ to the voltage
induced through~$e'$ by a voltage source between $e_-$ and $e_+$ was
used in \cite{NaddafSpencer}, Proposition~1, in
\cite{GloriaOttoVariance}, Lemma~2.4 and in \cite{BiskupSalviWolff12},
Proposition~2.5.
\end{rem}

The formula satisfied by the partial derivatives of $r\mapsto
i^{u,v}_r$ in Lemma~\ref{lem:regularite1} allows us to control
$i^{u,v}_r$ after a finite number of modifications of the individual
resistances.

\begin{lemm}
\label{lem:regularite2}
For any subset $S\subset\E_{1/2}$, define $r^{S\leftarrow r'}$ by:
\[
r^{S\leftarrow r'}(e)=\cases{
r'(e),  &
\quad$\mbox{if }  e\in S,$
\vspace*{2pt}\cr
r(e), & \quad$\mbox{else}$.}
\]
Then:
\begin{longlist}[(ii)]
\item[(i)] For any $e\in\E_{1/2}$, if $r'(e)\leq r(e)$,
\[
\bigl|i^{u,v}_{r}(e)\bigr|\leq\bigl|i^{u,v}_{r^{e\leftarrow r'}}(e)\bigr|
\leq\frac
{r(e)}{r'(e)}\bigl|i^{u,v}_{r}(e)\bigr|.
\]
\item[(ii)] Let $g(x,y)=\max\{\frac{x}{y},\frac{y}{x}\}$ for $x$ and
$y$ in $(0,+\infty)$. For any nonempty, finite subset $S\subset\E
_{1/2}$, any edge $e\in\E_{1/2}$, and any distinct vertices $u$ and
$v$, if $e\notin S$,
\[
\bigl\llvert i^{u,v}_{r^{S\leftarrow r'}}(e)\bigr\rrvert \leq\bigl\llvert
i^{u,v}_{r}(e)\bigr\rrvert + \biggl(\sum
_{e'\in S}\bigl\llvert i^{u,v}_{r}
\bigl(e'\bigr)\bigr\rrvert \biggr)\prod
_{e'\in S}g\bigl(r\bigl(e'\bigr),r'
\bigl(e'\bigr)\bigr),
\]
and if $e\in S$,
\[
\bigl\llvert i^{u,v}_{r^{S\leftarrow r'}}(e)\bigr\rrvert \leq \biggl(\sum
_{e'\in
S}\bigl\llvert i^{u,v}_{r}
\bigl(e'\bigr)\bigr\rrvert \biggr)\prod
_{e'\in S}g\bigl(r\bigl(e'\bigr),r'
\bigl(e'\bigr)\bigr).
\]
\end{longlist}
\end{lemm}

\begin{pf}
First, let us prove (i). Let $e\in\E_{1/2}$ and consider $\{
r(e'),e'\neq e\}$ fixed in $(0,\infty)^{\E_{1/2}}$. To simplify
notation, for $x>0$, define
\[
f(x):=i^{u,v}_{r^{e\leftarrow x}}(e)
\]
and
\[
g(x):=i^{e}_{r^{e\leftarrow x}}(e).
\]
Then one gets from Lemma~\ref{lem:regularite1},
\[
f'(x)=\frac{f(x)}{x}\bigl(g(x)-1\bigr).
\]
This is a homogeneous differential equation of order 1 on $(0,\infty)$
which implies that $f$ is of constant sign: either it is zero on
$(0,\infty)$, or it is positive on $(0,\infty)$, or it is negative on
$(0,\infty)$. Suppose that it is not identically zero and orient $e$ so
that $f$ is positive. Notice that $g(x)\in\,]0,1]$ for every $x>0$. Then
$f'$ is negative, which shows that for any $x_0\leq x_1$,
\[
f(x_1)\leq f(x_0).
\]
But also,
\[
(-\ln f)'(x)\leq\frac{1}{x}
\]
and thus
\[
f(x_0)\leq\frac{x_1}{x_0}f(x_1).
\]
This shows (i), and notably implies the following:
%
\begin{equation}
\label{eq:ee}
\bigl|i_{r^{e\leftarrow r'}}(e)\bigr|\leq
\bigl|i_{r}(e)\bigr|\max\biggl\{1,
\frac{r(e)}{r'(e)}\biggr\}.
\end{equation}
Now, let $e'$ and $e$ be distinct edges in $\E_{1/2}$. Using
Lemma~\ref{lem:regularite1}, and dropping the superscript $u,v$,
\[
\frac{\partial i_{r^{e'\leftarrow x}}(e)}{\partial x}= \frac
{i_{r^{e'\leftarrow x}}(e')}{x}i^{e'}_{r^{e'\leftarrow x}}(e).
\]
Recall from Lemma~\ref{lem:boundcurrent} that $|i^{e'}_{r^{e'\leftarrow
x}}(e)|$ is not larger than one, and from inequality~(\ref{eq:ee}),
\[
\bigl|i_{r^{e'\leftarrow x}}\bigl(e'\bigr)\bigr|\leq\bigl|i_{r}
\bigl(e'\bigr)\bigr|\max\biggl\{1,\frac{r(e')}{x}\biggr\}.
\]
Thus,
\begin{eqnarray*}
\bigl|i_{r^{e'\leftarrow r'}}(e)-i_{r}(e)\bigr|&\leq&\bigl|i_{r}
\bigl(e'\bigr)\bigr|\biggl\llvert \int_{r(e)}^{r(e')}
\frac{1}{x}\max\biggl\{1,\frac{r(e')}{x}\biggr\} \,dx\biggr\rrvert
\\
&\leq&\bigl|i_{r}\bigl(e'\bigr)\bigr|\max\bigl\{r
\bigl(e'\bigr),r'\bigl(e'\bigr)\bigr\}
\biggl\llvert \frac{1}{r(e')}-\frac
{1}{r'(e')}\biggr\rrvert
\\
&= &\bigl|i_{r}\bigl(e'\bigr)\bigr|\bigl(g\bigl(r
\bigl(e'\bigr),r'\bigl(e'\bigr)\bigr)-1
\bigr).
\end{eqnarray*}
Thus,
%
\begin{equation}
\label{eq:ee'}
\bigl|i_{r^{e'\leftarrow r'}}(e)\bigr|\leq \bigl|i_{r}(e)\bigr|+
\bigl|i_{r}
\bigl(e'\bigr)\bigr|\bigl(g\bigl(r\bigl(e'
\bigr),r'\bigl(e'\bigr)\bigr)-1\bigr).
\end{equation}
Notice that we have now established (ii) for sets $S$ of size 1.
Now, let us prove the first part of (ii) by induction on the size
of the set $S$. Let $e\notin S$. Using inequality~(\ref{eq:ee'}),
\begin{eqnarray*}
\bigl|i_{r^{S\leftarrow r'}}(e)\bigr|&=&\bigl\llvert i_{(r^{S\setminus\{e'\}\leftarrow
r'})^{e'\leftarrow r'}}(e)\bigr\rrvert
\\
&\leq&\bigl\llvert i_{r^{S\setminus\{e'\}\leftarrow r'}}(e)\bigr\rrvert +\bigl\llvert
i_{r^{S\setminus\{e'\}\leftarrow r'}}\bigl(e'\bigr)\bigr\rrvert \bigl(g\bigl(r
\bigl(e'\bigr),r'\bigl(e'\bigr)\bigr)-1
\bigr).
\end{eqnarray*}
From the induction hypothesis,
\[
\bigl\llvert i_{r^{S\setminus\{e'\}\leftarrow r'}}(e)\bigr\rrvert \leq \bigl|i_{r}(e)\bigr|+
\biggl(\sum_{e''\in S\setminus\{e'\}}\bigl\llvert i_{r}
\bigl(e''\bigr)\bigr\rrvert \biggr)\prod
_{e''\in S\setminus\{e'\}}g\bigl(r\bigl(e''
\bigr),r'\bigl(e''\bigr)\bigr)
\]
and
\[
\bigl\llvert i_{r^{S\setminus\{e'\}\leftarrow r'}}\bigl(e'\bigr)\bigr\rrvert
\leq\bigl|i_{r}\bigl(e'\bigr)\bigr| + \biggl(\sum
_{e''\in S\setminus\{e'\}}\bigl\llvert i_{r}\bigl(e''
\bigr)\bigr\rrvert \biggr)\prod_{e''\in S\setminus\{e'\}}g\bigl(r
\bigl(e''\bigr),r'\bigl(e''
\bigr)\bigr).
\]
Gathering terms, and noting that $g$ is not smaller than 1 allows to
complete the induction step. Finally, the second part of (ii) is
a consequence of the first part and inequality~(\ref{eq:ee'}). Indeed,
if $e\in S$,
\begin{eqnarray*}
\bigl|i_{r^{S\leftarrow r'}}(e)\bigr|&=& \bigl\llvert i_{(r^{S\setminus\{e\}\leftarrow
r'})^{e'\leftarrow r'}}(e)\bigr\rrvert
\\
&\leq&\bigl\llvert i_{r^{S\setminus\{e\}\leftarrow r'}}(e)\bigr\rrvert g\bigl(r(e),r'(e)
\bigr)
\\
&\leq& \biggl(\bigl|i_r(e)\bigr|+\sum_{e'\in S\setminus\{e\}}\bigl|i_r
\bigl(e'\bigr)\bigr|\prod_{e'\in
S\setminus\{e\}}g\bigl(r
\bigl(e'\bigr),r\bigl(e'\bigr)\bigr) \biggr)g
\bigl(r(e),r'(e)\bigr)
\\
&\leq& \sum_{e'\in S}\bigl|i_r
\bigl(e'\bigr)\bigr|\prod_{e'\in S}g\bigl(r
\bigl(e'\bigr),r\bigl(e'\bigr)\bigr).
\end{eqnarray*}
\upqed\end{pf}
%
\subsection{The random setting}
For any $e\in\E_{1/2}$, we let $\mu_e$ denote some probability measure
on $(0,\infty)$. The collection of resistances $r$ will be supposed to
be random with distribution $\PP:=\bigotimes_{e\in\E_{1/2}}\mu_e$.
\textit{Furthermore}, \textit{in the sequel},
$r'$ \textit{will usually denote an
independent copy of} $r$.

We shall always suppose that the resistances are square integrable.
Recall that the network is said to be \textit{elliptic} if there is a
constant $\Lambda>1$ such that $r\in[1,\Lambda]^{\E_{1/2}}$. This will
be a crucial assumption from Section~\ref{subsec:homogeneous} on.
Finally, we will use the notation:
\[
\forall e\in\E, \qquad m_p(e)=\EE\bigl[\bigl|r(e)-\EE\bigl(r(e)
\bigr)\bigr|^p\bigr]^{1/p}.
\]

\section{The Walsh decomposition}
\label{sec:ES}
\subsection{Definition and basic properties}
\label{sec:ESbasics}
For any $e\in\E_{1/2}$ let $\Delta_e$ be the following operator on
$L^2(\RR^{\E_{1/2}},\PP)$:
\[
\Delta_ef(r)=f(r)-\int f(r) \,d\mu_e\bigl(r(e)\bigr).
\]
From now on, $S\subset\E_{1/2}$ will always mean that $S$ is a \textit{finite} subset of $\E_{1/2}$. For $S\subset\E_{1/2}$, we shall denote
by $r_S$ the collection $(r(e))_{e\in S}$ of random variables (which is
empty if $S$ is empty). Let $f$ be in $L^2(\RR^{\E_{1/2}},\PP)$ and
notice that
\[
\EE\bigl[f(r)|r_S\bigr]=\int f(r) \bigotimes
_{e\in S^c}\,d\mu_e\bigl(r(e)\bigr).
\]
Then, for any $S\subset\E_{1/2}$ we define
\[
f_S(r)=\sum_{T\subset S}(-1)^{|S\setminus T|}
\EE\bigl[f(r)|r_T\bigr].
\]
Notice that $f_{\varnothing}=\EE(f)$. It is easy to see that an
alternative definition is
\[
f_S(r)=\EE\biggl[\biggl(\prod_{e\in S}
\Delta_e\biggr) f(r)\Big|r_S\biggr],
\]
with the usual convention that when $S$ is empty, the product of
operators over~$S$ is the identity. Then $(f_S)_{S\subset\E_{1/2}}$ is
an orthogonal decomposition of $f$ known as the Efron--Stein or the
(generalized) Walsh decomposition (cf. \cite{Bourgain79,EfronStein81,Hatami12} and \cite{Mossel10},  e.g.). The
basic properties of this decomposition are gathered in the following
proposition, where infinite sums in $L^2(\RR^{\E_{1/2}},\PP)$ are
understood as follows: $\sum_Sf_S$ is the limit in $L^2(\RR^{\E
_{1/2}},\PP)$ of the net $S\mapsto\sum_{T\subset S}f_T$, defined on
the set of finite subsets of $\E_{1/2}$ with inclusion as partial order
(in other words, this corresponds to unconditional summability in $L^2$).

\begin{prop}
\label{prop:ES}
For any $f$ and $g$ in $L^2(\RR^{\E_{1/2}},\PP)$,
\begin{eqnarray*}
f &=& \sum_Sf_S,
\\
\EE(fg) &=& \sum_S\EE(f_Sg_S),
\end{eqnarray*}
and thus
\[
S\neq S'\quad \Rightarrow \quad \EE(f_Sg_{S'})=0.
\]
Furthermore, for any $e\in\E_{1/2}$,
\[
\Delta_ef=\sum_{S\ni e}f_S.
\]
As a consequence, for any integer $k\geq1$,
\[
\sum_{{e_1,\ldots, e_k}\in\E_{1/2}}\Biggl\llVert \Biggl(\prod
_{i=1}^k\Delta _{e_i} \Biggr)f\Biggr
\rrVert ^2=\sum_S|S|^k
\|f_S\|_2^2.
\]
\end{prop}

\begin{pf}
Let us first suppose that $\E_{1/2}$ is finite. Then Proposition~\ref
{prop:ES} is well known (cf. \cite{Bourgain79}) but we shall quickly
recall the proof for the sake of completeness.

For a subset $S$ any subset of $\E_{1/2}$, let $L_S$ be the operator on
$L^2(\RR^{\E_{1/2}})$ defined by
%
\begin{equation}
\label{eq:defLS}
L_Sf(r)=\int f(r) \prod
_{e\in S}\,d\mu_e\bigl(r(e)\bigr).
\end{equation}
Let $1$ denote the identity operator. Notice that $L_{\{e\}}$ and $L_{\{
e'\}}$ commute for any $e$ and $e'$. Since $\Delta_{e'}=1-L_{\{e\}}$,
$\Delta_{e'}$ and $L_{\{e\}}$ commute and
\[
1=\prod_{e\in\E_{1/2}}(\Delta_{\{e\}}+L_{\{e\}})=
\sum_{S\subset\E
_{1/2}}L_{S^c}\prod
_{e\in S}\Delta_{\{e\}}.
\]
Since
\[
f_S= \biggl(\prod_{e\in S}
\Delta_{\{e\}} \biggr)f,
\]
this shows that $f_\varnothing=\EE(f)$ and $f=\sum_{S\subset\E
_{1/2}}f_S$. Now, remark that for any edge $e$,
\[
L_{\{e\}}\Delta_e=0.
\]
Thus, for any $S$ and any $e\in S$,
\[
L_{\{e\}}f_S=0.
\]
This implies that $\Delta_ef=\sum_{S\ni e}f_S$. Now, if $S\neq S'$,
suppose, for instance, that there is some $e\in S\setminus S'$:
\[
\EE[f_Sg_{S'}]=L_{\E_{1/2}}(f_Sg_{S'})=L_{\E_{1/2}}L_{\{e\}
}(f_Sg_{S'})=L_{\E_{1/2}}
\bigl(g_{S'}L_{\{e\}}(f_S)\bigr)=0.
\]
This implies
\[
\EE(fg)=\sum_S\EE(f_Sg_S).
\]
Finally,
\begin{eqnarray*}
\sum_{{e_1,\ldots, e_k}\in\E_{1/2}}\Biggl\llVert \Biggl(\prod
_{i=1}^k\Delta _{e_i} \Biggr)f\Biggr
\rrVert ^2&=&\sum_{{e_1,\ldots, e_k}\in\E_{1/2}}\biggl\llVert \sum
_{S\supset\{e_1,\ldots,e_k\}}f_S\biggr\rrVert ^2
\\
&=&\sum_{{e_1,\ldots, e_k}\in\E_{1/2}} \sum_{S\supset\{e_1,\ldots,e_k\}
}
\|f_S\|^2
\\
&=&\sum_{S\subset\E_{1/2}}\sum_{\{e_1,\ldots,e_k\}\subset S}
\|f_S\|^2
\\
&=&\sum_{S\subset\E_{1/2}}|S|^k\|f_S
\|^2.
\end{eqnarray*}

Now, let us suppose that $\E_{1/2}$ is countable, and take some
exhaustion $(E_n)_{n\geq0}$ of the edges: $E_n$ is finite for any $n$,
$E_n \subset E_{n+1}$ and $\E_{1/2}=\bigcup_{n}E_n$. Denote by $f_n$ the
conditional expectation of $f$ with respect to $r_{E_n}$. We have
\[
f_n=L_{E_n^c}f
\]
and thus,
%
\begin{equation}
\label{eq:fnS}
(f_n)_S=f_S
\II_{S\subset E_n},
\end{equation}
which implies
\[
f_n=\sum_{S\subset E_n}f_S \quad\mbox{and}\quad\|f_n\|_2^2=\sum
_{S\subset E_n}\| f_S\|_2^2.
\]
Since $(f_n)_{n\in\NN}$ converges to $f$ in $L^2$, we know that $\|f_n\|
_2^2$ converges to $\|f\|_2^2$. It is then standard to see that $\sum_Sf_S$ forms a Cauchy net and that its limit is the same as the limit
of $f_n$, that is $f$. It shows notably that $\EE(f^2)=\sum_S\EE
(f_S^2)$, from which one derives $\EE(fg)=\sum_S\EE(f_Sg_S)$. All the
other properties can then easily be derived by standard limit arguments.
\end{pf}

A trivial consequence of Proposition~\ref{prop:ES} is the Efron--Stein
inequality (cf. \cite{EfronStein81}).

\begin{coro}[(Efron--Stein's inequality)]
\label{coro:ESI}
For any $f$ in $L^2(\RR^{\E_{1/2}},\PP)$
\[
\Var(f)\leq\sum_{e\in\E_{1/2}}\|\Delta_ef
\|^2.
\]
\end{coro}

\begin{pf}
Since the Efron--Stein decomposition is orthogonal and $f_\varnothing=\EE(f)$,
\begin{eqnarray*}
\Var(f)&=&\EE\bigl(f^2\bigr)-\EE(f)^2
\\
&=&\sum_{S}\EE\bigl(f_S^2
\bigr)-f_\varnothing^2
\\
&=&\sum_{S\neq\varnothing}\EE\bigl(f_S^2
\bigr).
\end{eqnarray*}
On the other hand, since $\Delta_ef=\sum_{S\ni e}f_S$,
\begin{eqnarray*}
\sum_{e\in\E_{1/2}}\|\Delta_ef
\|^2&=&\sum_{e\in\E_{1/2}}\sum
_{S\ni
e}\EE\bigl(f_S^2\bigr)
\\
&=&\sum_{S}\sum_{e\in S}
\EE\bigl(f_S^2\bigr)
\\
&=&\sum_{S}|S|\EE\bigl(f_S^2
\bigr)
\\
&\geq&\sum_{S\neq\varnothing}\EE\bigl(f_S^2
\bigr).
\end{eqnarray*}
\upqed\end{pf}
It is also clear that the Efron--Stein inequality is an equality if and
only if $f=\sum_{|S|\leq1}f_S$. This\vspace*{1pt} means that $f$ is a constant plus
a sum of independent random variables. One sees also that if the
variance of $f$ is concentrated on functions $f_S$ such that $S$ is
small, then the Efron--Stein inequality is sharp up to a multiplicative
constant.

There are a number of models of statistical physics flavour where the
Efron--Stein inequality is not sharp. Chatterjee (cf. \cite{ChatterjeeChaos})
calls this phenomenon ``superconcentration.'' This
holds, for instance, for the first passage percolation time between two
distant points on $\ZZ^d$ when $d\geq2$ (cf. \cite{BenjaminiKalaiSchramm99,BenaimRossignol08}), which may be defined in
our setting as
%
\begin{equation}
\label{eq:Tr}
T_r(u,v):=\inf_{\gamma \dvtx  u\to v}\sum
_{e\in\gamma} r(e),
\end{equation}
where the infimum is over all paths from $u$ to $v$. Since
super-concentration implies that some part of the variance of $f$ is
concentrated on large sets, one sees that, informally, it is related to
high complexity (or high nonlinearity) of the function $f$.

It is also related to some noise-sensitivity of the function. Indeed,
there is a close link between the Walsh decomposition and a notion of
noise introduced by~\cite{BenjaminiKalaiSchramm99}. Let $r$ and $r'$ be
two independent random variables with the same distribution $\PP
:=\bigotimes_{e\in E_{1/2}}\mu_e$. Let $\eps\in\,]0,1[\,$.\vspace*{1pt}
One constructs a
noisy version $r^\eps$ of $r$ by replacing with probability $\eps$, at
random and independently for any edge $e$, the variable $r(e)$ by its
independent copy $r'(e)$.

\begin{prop}
\label{prop:noise}
For any $f\in L^2(\RR^{\E_{1/2}},\PP)$,
\[
\EE\bigl[f\bigl(r^\eps\bigr)|r\bigr]=\sum
_S(1-\eps)^{|S|}f_S(r),
\]
and thus,
\[
\Cov\bigl(f\bigl(r^\eps\bigr),f(r)\bigr)=\sum
_{S\neq\varnothing}(1-\eps)^{|S|}\|f_S
\|_2^2.
\]
\end{prop}

\begin{pf}
Let $r$ and $r'$ be two independent copies of law $\PP$. Let $S_\eps$
be the (possibly infinite) random\vspace*{1pt} subset of $\E$ drawn at random as
follows: $(\II_{e\in S_{\eps}})_{e\in\E_{1/2}}$ are i.i.d. with
distribution Bernoulli of parameter $\eps\in[0,1]$, independent of
$(r,r')$. Now, we define the following linear operators from $L^2(\RR
^{\E_{1/2}})$ to $L^2(\RR^{\E_{1/2}}\times\RR^{\E_{1/2}})$:
\[
L^{r'}_{\{e\}}f(r)=f\bigl(r^{e\leftarrow r'}\bigr).
\]
Then the noisy version of $f(r)$ may be written as
\[
f\bigl(r^\eps\bigr)=\prod_{e\in S_\eps}L^{r'}_{\{e\}}f(r).
\]
Notice that
\[
\cases{\ds \forall e\neq e', &
$L^{r'}_{\{e\}}L_{\{e'\}}=L_{\{e'\}}L^{r'}_{\{
e\}}$,
\vspace*{2pt}\cr
\ds\forall e, & $L^{r'}_{\{e\}}L_{\{e\}}=L_{\{e\}}$,\vspace*{2pt}
\cr
\ds\forall e,  & $L_{\{e\}}L^{r'}_{\{e\}}=L^{r'}_{\{e\}}$.}
\]
Thus,
\[
L^{r'}_{\{e\}}\Delta_e=L^{r'}_{\{e\}}-L_{\{e\}}.
\]
Whence
\begin{eqnarray*}
f\bigl(r^\eps\bigr)&=&\sum_{S\subset\E_{1/2}}\prod
_{e\in S_\eps}L^{r'}_{\{e\}
}f_S(r)
\\
&=&\sum_{S\subset\E_{1/2}}\prod_{e\in S_\eps}L^{r'}_{\{e\}}
\prod_{e\in
S^c}L_{\{e\}}\prod
_{e\in S}\Delta_e f(r)
\\
&=&\sum_{S\subset\E_{1/2}}\prod_{e\in S^c}L_{\{e\}}
\prod_{e\in S\cap
S_\eps}\bigl(L^{r'}_{\{e\}}-L_{\{e\}}
\bigr)\prod_{e\in S\setminus S_\eps}\Delta_e f(r).
\end{eqnarray*}
Now, denote by $L'_{\{e\}}$ the operator on $L^2(\RR^{\E_{1/2}}\times\RR
^{\E_{1/2}})$ which integrates $r'(e)$. Notice that
\[
L'_{\{e\}}\bigl(L^{r'}_{\{e\}}-L_{\{e\}}
\bigr)=0 \quad \mbox{and} \quad L'_{\{e\}}\Delta _e=
\Delta_e.
\]
Then
\begin{eqnarray*}
\EE\bigl[f\bigl(r^\eps\bigr)|r,S_\eps\bigr]&=&\prod
_{e\in\E_{1/2}}L'_{\{e\}}\bigl(
\bigl(r,r'\bigr)\mapsto f\bigl(r^\eps\bigr)\bigr)
\\
&=&\sum_{S\subset\E_{1/2}}\prod_{e\in S^c}L_{\{e\}}
\prod_{e\in
S\setminus S_\eps}\Delta_e f(r)
\II_{S\cap S_\eps=\varnothing}
\\
&=&\sum_{S\subset\E_{1/2}}\prod_{e\in S^c}L_{\{e\}}
\prod_{e\in S}\Delta _e f(r)
\II_{S\cap S_\eps=\varnothing}.
\end{eqnarray*}
Thus,
\[
\EE\bigl[f\bigl(r^\eps\bigr)|r\bigr]=\sum
_{S\subset\E_{1/2}} f_S(r)\PP(S\cap S_\eps =
\varnothing)=\sum_{S\subset\E_{1/2}}(1-\eps)^{|S|}f_S(r).\quad\qed
\]
\noqed\end{pf}

To see another, closely related, interpretation of sensitivity to
noise, called chaos, see \cite{ChatterjeeChaos}. Contrarily to what
happens to the first passage percolation times, simulations suggest
that the minimal current is extremely immune to noise. This tends to
suggest that the Efron--Stein inequality could always be sharp in this
context, and this is what we shall prove in the following sections. We
shall in fact prove much more in the context of \textit{current-homogeneous graphs}, namely that the Walsh decomposition is
concentrated not only on sets of small size, but already on sets of
small diameter.

Finally, to end the parallel between first passage percolation and
effective resistance, note that there is a means to interpolate between
those two quantities. Indeed, let us define, for $p\in[1,2]$:
\[
\R^{u,v}_r(p)=\inf_{\theta\dvtx u\to v}\sum
_{e\in\E_{1/2}} r(e)\bigl|\theta (e)\bigr|^p,
\]
where the infimum is over all unit flows from $u$ to $v$, and recall
the definition~(\ref{eq:Tr}) of $T_r(u,v)$, the minimum passage time
from $u$ to $v$. Then $\R^{u,v}_r(1)=T_r(u,v)$ and $\R^{u,v}_r(2)=\R
^{u,v}(r)$. The quantity $\R^{u,v}_r(p)$ is called the $p$-resistance
between $u$~and~$v$. Since the distribution of the Walsh decomposition
is dramatically different when $p$ equals 1 or 2, it would be
interesting to investigate the evolution of the Walsh decomposition of
$\R^{u,v}_r(p)$ when $p$ varies continuously from $2$ to $1$.

\subsection{Concentration of the Walsh decomposition on low levels}
\label{sec:lowlevels}
First, let us study the bound given by the Efron--Stein inequality.

\begin{lemm}
\label{lem:thetae4Deltae}
For any $e\in\E_{1/2}$,
\[
\alpha_-(e)\EE\bigl[r^2(e) \bigl(i_{r}^{u,v}(e)
\bigr)^4\bigr]\leq\bigl\|\Delta_e\R^{u,v}\bigr\|
_2^2\leq\alpha_+(e)\EE\bigl[r^2(e)
\bigl(i_{r}^{u,v}(e)\bigr)^4\bigr],
\]
where
\[
\alpha_-(e)=\frac{\EE [(r(e)-r'(e))_+^2\min \{1,{1}/({r^4(e)}) \} ]}{\EE [\max \{r^2(e),{1}/({r^2(e)}) \} ]}
\]
and
\[
\alpha_+(e)=\frac{\EE [(r'(e)-r(e))_+^2\max \{1,
{1}/({r^4(e)}) \} ]}{\EE [\min\{r^2(e),{1}/({r^2(e)})\}
 ]}.
\]
\end{lemm}

\begin{pf}
Let $r$ and $r'$ be two independent random variables with the same
distribution $\PP:=\bigotimes_{e\in E_{1/2}}\mu_e$. Inequality
(\ref{eq:majoDeltaeR2sided}) gives
\begin{eqnarray*}
\bigl\|\Delta_e\R^{u,v}\bigr\|_2^2&=&\EE
\bigl[\bigl(\R^{u,v}(r)-\R^{u,v}\bigl(r^{e\leftarrow
r'}\bigr)
\bigr)_+^2\bigr]
\\
&\geq&\EE\bigl[\bigl(r(e)-r'(e)\bigr)_+^2i_r^{u,v}(e)^4
\bigr].
\end{eqnarray*}
Now, Lemma~\ref{lem:regularite2} allows to decouple positive functions
of $r(e)$ and powers of $|i_r(e)|$. Let $F$ and $G$ be nonnegative
functions on $(0,+\infty)$ and $p$ be a positive real number. Then
%
\begin{eqnarray}
&& \EE\bigl[F\bigl(r(e)\bigr)\bigl|i_r^{u,v}(e)\bigr|^p
\bigr]
\nonumber
\\[-8pt]
\label{eq:decoupling}
\\[-8pt]
\nonumber
&&\qquad\leq\EE\bigl[G\bigl(r(e)\bigr)\bigl|i_r^{u,v}(e)\bigr|^p
\bigr] \frac{\EE
[F(r(e))\max\{1,{1}/({r^p(e)})\}]}{\EE[G(r(e))\min\{1,
{1}/({r^p(e)})\}]}.
\end{eqnarray}
Indeed, using Lemma~\ref{lem:regularite2}, for any $r$,
\[
\bigl|i_r^{u,v}(e)\bigr|\leq\max\biggl\{1,\frac{1}{r(e)}\biggr
\}\bigl|i_{r^{e\leftarrow
1}}^{u,v}(e)\bigr|.
\]
Thus,
\begin{eqnarray*}
\EE\bigl[F\bigl(r(e)\bigr)\bigl|i_r^{u,v}(e)\bigr|^p
\bigr]&\leq&\EE\biggl[F\bigl(r(e)\bigr)\max\biggl\{1,\frac{1}{r^p(e)}\biggr\}
\bigl|i_{r^{e\leftarrow1}}^{u,v}(e)\bigr|^p\biggr]
\\
&=&\EE\biggl[F\bigl(r(e)\bigr)\max\biggl\{1,\frac{1}{r^p(e)}\biggr\}\biggr]\EE
\bigl[\bigl|i_{r^{e\leftarrow
1}}^{u,v}(e)\bigr|^p\bigr],
\end{eqnarray*}
since $r(e)$ and $i_{r^{e\leftarrow1}}^{u,v}(e)$ are independent. Similarly,
\begin{eqnarray*}
\EE\bigl[G\bigl(r(e)\bigr)\bigl|i_r^{u,v}(e)\bigr|^p
\bigr]&\geq& \EE\biggl[G\bigl(r(e)\bigr)\min\biggl\{1,\frac
{1}{r^p(e)}\biggr
\}\bigl|i_{r^{e\leftarrow1}}^{u,v}(e)\bigr|^p\biggr]
\\
&=&\EE\biggl[G\bigl(r(e)\bigr)\min\biggl\{1,\frac{1}{r^p(e)}\biggr\}\biggr]\EE
\bigl[\bigl|i_{r^{e\leftarrow
1}}^{u,v}(e)\bigr|^p\bigr].
\end{eqnarray*}
%
Thus,
\[
\bigl\|\Delta_e\R^{u,v}\bigr\|_2^2 \geq
\frac{\EE [(r(e)-r'(e))_+^2\min \{
1,{1}/({r^4(e)}) \} ]}{\EE [\max \{r^2(e),{1}/({r^2(e)}) \} ]}\EE\bigl[r^2(e)i_{r^{e\leftarrow1}}^{u,v}(e)^4
\bigr].
\]
On the other hand,
\begin{eqnarray*}
\bigl\|\Delta_e\R^{u,v}\bigr\|_2^2&=&\EE
\bigl[\bigl(\R^{u,v}\bigl(r^{e\leftarrow r'}\bigr)-\R ^{u,v}(r)
\bigr)_+^2\bigr]
\\
&\leq&\EE\bigl[\bigl(r'(e)-r(e)\bigr)_+^2i_r^{u,v}(e)^4
\bigr]
\\
&\leq&\frac{\EE [(r'(e)-r(e))_+^2\max \{1,{1}/({r^4(e)}) \} ]}{\EE [\min\{r^2(e),{1}/({r^2(e)})\}
 ]}\EE
\bigl[r^2(e)i_{r}^{u,v}(e)^4\bigr].
\end{eqnarray*}
\upqed\end{pf}
The following theorem shows that the Walsh decompositions of
point-to-point effective resistances are uniformly concentrated (in
terms of the $L^2$-norm) on sets of small size.

\begin{theo}
\label{theo:EStight}
There is a universal constant $C\in(0,+\infty)$ such that if one
defines in $[0,+\infty]$:
\[
K(\mu)=C\sup_e\bigl(\EE\bigl[r^8(e)\bigr]+\EE
\bigl[r^{-8}(e)\bigr]\bigr)^6\sup_e
\frac{\EE
[(r(e)-r'(e))^2 ({1}/({r^6(e)})+r^6(e) ) ]}{\EE
[(r(e)-r'(e))_+^2\min \{1,{1}/({r^4(e)}) \} ]},
\]
then for any graph $\G$ and any pair of vertices $(u,v)$,
\[
\sum_{S}|S|^2\bigl\|\R^{u,v}_S
\bigr\|_2^2\leq K(\mu)\sum_{S\neq\varnothing}\bigl\|
\R ^{u,v}_S\bigr\|_2^2.
\]
Consequently, for any $k\geq1$, any graph $\G$ and any pair of
vertices $(u,v)$
\[
\sum_{|S|\geq k}\bigl\|\R^{u,v}_S
\bigr\|_2^2\leq\frac{K(\mu)}{k^2}\sum
_{S\neq\varnothing}\bigl\|\R^{u,v}_S\bigr\|_2^2.
\]
\end{theo}

\begin{pf}
We fix $u$ and $v$ and drop the superscript $u,v$. Proposition~\ref{prop:ES} implies
\[
\sum_{S}|S|^2\|\R_S
\|_2^2=\sum_{e,e'}\|
\Delta_e\Delta_{e'}\R\|_2^2.
\]
The first step is to prove the following:
%
\begin{equation}
\label{eq:majoH2H1}
\sum_{e\neq e'}\|\Delta_e
\Delta_{e'}\R\|_2^2\leq D(\mu)\sum
_{e}\bigl\|\Delta _e\R^{u,v}
\bigr\|_2^2.
\end{equation}
Suppose for the moment that (\ref{eq:majoH2H1}) is true. Then
%
\begin{eqnarray}
\nonumber
\sum_{S}|S|^2 \cdot \bigl\|
\R^{u,v}_S\bigr\|_2^2&=&\sum
_{e,e'}\bigl\|\Delta_e\Delta _{e'}
\R^{u,v}\bigr\|_2^2
\\
\nonumber
&=&\sum_{e\neq e'}\bigl\|\Delta_e
\Delta_{e'}\R^{u,v}\bigr\|_2^2+\sum
_{e}\bigl\|\Delta_e\R^{u,v}
\bigr\|_2^2
\\[-8pt]
\label{eq:card2card1}
\\[-8pt]
\nonumber
&\leq&\bigl(D(\mu)+1\bigr)\sum_{e}\bigl\|
\Delta_e\R^{u,v}\bigr\|_2^2
\\
\nonumber
&=&\bigl(D(\mu)+1\bigr)\sum_{S}|S|
\bigl\|\R^{u,v}_S\bigr\|_2^2.
\end{eqnarray}
Then, for any $k\in\NN^*$,
\begin{eqnarray*}
\sum_{|S|\geq k}|S|\bigl\|\R^{u,v}_S
\bigr\|_2^2&\leq&\frac{1}{k}\sum
_{S}|S|^2\bigl\|\R ^{u,v}_S
\bigr\|_2^2
\\
&\leq&\frac{(D(\mu)+1)}{k}\sum_{S}|S|\bigl\|
\R^{u,v}_S\bigr\|_2^2.
\end{eqnarray*}
Thus, for every\vspace*{-3pt} $k\geq2(D(\mu)+1)$,
\[
\sum_{S}|S|\bigl\|\R^{u,v}_S
\bigr\|_2^2\leq k\sum_{0<|S|<k}\bigl\|
\R^{u,v}_S\bigr\| _2^2+
\frac{1}{2} \sum_{S}|S|\bigl\|
\R^{u,v}_S\bigr\|_2^2,
\]
whence,\vspace*{-2pt} for every $k\geq2(D(\mu)+1)$,
\[
\sum_{S}|S|\bigl\|\R^{u,v}_S
\bigr\|_2^2\leq2k\sum_{0<|S|<k}\bigl\|
\R^{u,v}_S\bigr\|_2^2.
\]
Notably,\vspace*{-4pt}
\[
\sum_{S}|S|\bigl\|\R^{u,v}_S
\bigr\|_2^2\leq2\bigl\lceil2\bigl(D(\mu)+1\bigr)\bigr\rceil\sum
_{S\neq\varnothing}\bigl\|\R^{u,v}_S
\bigr\|_2^2.
\]
Plugging this into inequality (\ref{eq:card2card1}) ends the proof of
the first inequality of the theorem\vspace*{-4pt} with
\[
K(\mu)= 2\bigl\lceil2\bigl(D(\mu)+1\bigr)\bigr\rceil.
\]

It remains to prove (\ref{eq:majoH2H1}). Let us present first the main
idea in the elliptic setting. We know from Lemma~\ref{lem:regularite1} that
%
\begin{equation}
\label{eq:dertheta2}
\partial^2_{e,e'}\R(r)=\partial_e
\bigl[i_r\bigl(e'\bigr)\bigr]^2=2i_r
\bigl(e'\bigr)i_r(e)\frac
{i_r^{e}(e')}{r(e)}.
\end{equation}
Approximating $\Delta_e\Delta_{e'}\R^{u,v}$ by $\partial^2_{e,e'}\R
(r)$, one gets
\begin{eqnarray*}
\sum_{e\neq e'}\bigl\|\Delta_e
\Delta_{e'}\R^{u,v}\bigr\|_2^2&\lesssim&
\sum_{e\neq e'}\EE\bigl[i_r
\bigl(e'\bigr)^2i_r(e)^2
\bigl(i_r^{e}\bigl(e'\bigr)
\bigr)^2\bigr]
\\[-2pt]
&\lesssim& \sum_{e\neq e'}\EE\bigl[\bigl(i_r
\bigl(e'\bigr)^4+i_r(e)^4
\bigr) \bigl(i_r^{e}\bigl(e'\bigr)
\bigr)^2\bigr].
\end{eqnarray*}
The reciprocity law~(\ref{eq:reciprocity}) gives that $i_r^{e}(e')$ and
$i_r^{e'}(e)$ are of the same\vspace*{-2pt} order:
\begin{eqnarray*}
\sum_{e\neq e'}\bigl\|\Delta_e
\Delta_{e'}\R^{u,v}\bigr\|_2^2&\lesssim&
\sum_{e\neq e'}\EE\bigl[i_r(e)^4
\bigl(i_r^{e}\bigl(e'\bigr)
\bigr)^2\bigr]+\sum_{e\neq e'}\EE
\bigl[i_r\bigl(e'\bigr)^4
\bigl(i_r^{e'}(e)\bigr)^2\bigr]
\\[-3pt]
&= &2\sum_{e\neq e'}\EE\bigl[i_r(e)^4
\bigl(i_r^{e}\bigl(e'\bigr)
\bigr)^2\bigr]
\\[-3pt]
&\lesssim&\sum_{e}\EE\biggl[i_r(e)^4
\sum_{e'}r\bigl(e'\bigr)
\bigl(i_r^{e}\bigl(e'\bigr)
\bigr)^2\biggr],
\end{eqnarray*}
but $\sum_{e'}r(e')(i_r^{e}(e'))^2$ is the effective resistance from
$e_-$ to $e_+$, which is of order at most 1 [in fact at most $r(e)$]. Thus,\vspace*{-2pt}
\begin{eqnarray*}
\sum_{e\neq e'}\bigl\|\Delta_e
\Delta_{e'}\R^{u,v}\bigr\|_2^2&\lesssim&
\sum_{e}\EE \bigl[i_r(e)^4
\bigr]
\\[-3pt]
&\lesssim&\sum_{e}\bigl\|\Delta_e
\R^{u,v}\bigr\|_2^2,
\end{eqnarray*}
from\vadjust{\goodbreak} Lemma~\ref{lem:thetae4Deltae}.

Now, let us enter the details of the proof of (\ref{eq:majoH2H1}) in
the general case. Let $r$ and $r'$ be two independent random variables
with the same distribution $\PP:=\bigotimes_{e\in E_{1/2}}\mu_e$.
Remark that for $e\neq e'$:
\begin{eqnarray*}
&&\Delta_e\Delta_{e'}\R(r)
\\
&& \qquad =\EE \bigl[\R(r)-\R\bigl(r^{e'\leftarrow r'}\bigr)-\R\bigl(r^{e\leftarrow
r'}
\bigr)+\R\bigl(r^{\{e,e'\}\leftarrow r'}\bigr)| r \bigr]
\\
&&\qquad =\EE \biggl[\int_{r'(e)}^{r(e)}\!\int
_{r'(e')}^{r(e')}\partial ^2_{e,e'}\R
\bigl(r^{(e,e')\leftarrow(x,y)}\bigr) \,dx\,dy\Big| r \biggr].
\end{eqnarray*}
Thus,
\begin{eqnarray*}
&&\|\Delta_e\Delta_{e'}\R\|_2^2
\\
&&\qquad \leq \EE \Bigl[\bigl(r(e)-r'(e)\bigr)^2\bigl(r
\bigl(e'\bigr)-r'\bigl(e'\bigr)
\bigr)^2 \mathop{\sup_{x\in[r(e),r'(e)]}}_{y\in[r(e'),r'(e')]} \bigl(
\partial ^2_{e,e'}\R\bigl(r^{(e,e')\leftarrow(x,y)}\bigr)
\bigr)^2 \Bigr],
\end{eqnarray*}
where we make the abuse of notation of writing $[a,b]$ for $[\min\{a,b\}
,\max\{a,b\}]$.
Lemma~\ref{lem:regularite2} shows that for $x$ in $[r(e),r'(e)]$ and
$y$ in $[r(e'),r'(e')]$,
\[
\bigl|i_{r^{(e,e')\leftarrow(x,y)}}(e)\bigr|\leq \bigl[\bigl|i_{r^{(e,e')\leftarrow1}}(e)\bigr|
+\bigl|i_{r^{(e,e')\leftarrow1}}
\bigl(e'\bigr)\bigr| \bigr]g(1,x)g(1,y),
\]
and the same bound holds for $|i_{r^{(e,e')\leftarrow(x,y)}}(e)|$.
Furthermore, using
Lemma~\ref{lem:regularite1} and the reciprocity law~(\ref{eq:reciprocity}),
\begin{eqnarray*}
\frac{1}{x}\bigl|i_{r^{(e,e')\leftarrow(x,y)}}^{e}\bigl(e'
\bigr)\bigr|&\leq&\frac
{1}{x}\bigl|i_{r^{(e,e')\leftarrow(x,1)}}^{e}
\bigl(e'\bigr)\bigr|\max\biggl\{1,\frac{1}{y}\biggr\}
\\
&=&\bigl|i_{r^{(e,e')\leftarrow(x,1)}}^{e'}(e)\bigr|\max\biggl\{1,\frac{1}{y}\biggr\}
\\
&\leq& \bigl|i_{r^{(e,e')\leftarrow1}}^{e'}(e)\bigr|\max\biggl\{1,\frac{1}{x}
\biggr\}\max\biggl\{ 1,\frac{1}{y}\biggr\}
\\
&=&\bigl|i_{r^{(e,e')\leftarrow1}}^{e}\bigl(e'\bigr)\bigr|\max\biggl\{1,
\frac{1}{x}\biggr\}\max\biggl\{ 1,\frac{1}{y}\biggr\}.
\end{eqnarray*}
Thus, using $(a+b)^4\leq8(a^4+b^4)$,
\[
\|\Delta_e\Delta_{e'}\R\|_2^2
\leq2A(e)A\bigl(e'\bigr)\EE \bigl[\bigl(i_{r^{(e,e')\leftarrow1}}(e)^4+i_{r^{(e,e')\leftarrow
1}}
\bigl(e'\bigr)^4\bigr)i_{r^{(e,e')\leftarrow1}}^{e}
\bigl(e'\bigr)^2 \bigr],
\]
where
\[
A(e)= 4\EE \biggl[\bigl(r(e)-r'(e)\bigr)^2 \biggl(
\frac{1}{r^6(e)}+r^6(e) \biggr) \biggr].
\]
Then we use again the same decoupling device based on Lemma~\ref
{lem:regularite2}, in the spirit of (\ref{eq:decoupling}). We get
\[
\bigl|i_{r^{(e,e')\leftarrow1}}^{e}\bigl(e'\bigr)\bigr|\leq
\frac{|i_r^e(e')|}{r(e)}\max\bigl\{ 1,r(e)\bigr\}\max\bigl\{1,r\bigl(e'
\bigr)\bigr\}
\]
and
\[
\bigl|i_{r^{(e,e')\leftarrow1}}(e)\bigr|\leq \bigl[\bigl|i_r(e)\bigr|+\bigl|i_r
\bigl(e'\bigr)\bigr| \bigr]g\bigl(r(e),1\bigr)g\bigl(r
\bigl(e'\bigr),1\bigr).
\]
Thus,
\begin{eqnarray*}
&&\EE \bigl[\bigl(i_{r^{(e,e')\leftarrow1}}(e)^4+i_{r^{(e,e')\leftarrow
1}}
\bigl(e'\bigr)^4\bigr)i_{r^{(e,e')\leftarrow1}}^{e}
\bigl(e'\bigr)^2 \bigr]
\\
&&\qquad \leq B(e)B\bigl(e'\bigr)\EE\biggl[\bigl(i_r^4(e)+i_r^4
\bigl(e'\bigr)\bigr)\frac
{i_r^e(e')^2}{r(e)^2}r(e)r\bigl(e'
\bigr)\biggr],
\end{eqnarray*}
where
\[
B(e)=\frac{4}{\EE\min\{{1}/({r^7(e)}),r^3(e)\}}.
\]
Whence
%
\begin{equation}
\label{eq:MajoDeltaee'}\qquad \|\Delta_e\Delta_{e'}\R
\|_2^2\leq2A(e)A\bigl(e'\bigr)B(e)B
\bigl(e'\bigr)\EE \biggl[\bigl(i_r^4(e)+i_r^4
\bigl(e'\bigr)\bigr)\frac{i_r^e(e')^2}{r(e)^2}r(e)r\bigl(e'
\bigr)\biggr].
\end{equation}
Now notice that for $e$ fixed,
\[
\sum_{e'}i_r^e
\bigl(e'\bigr)^2r\bigl(e'\bigr)=
\R^e(r)\leq r(e).
\]
We obtain
\begin{eqnarray*}
\sum_{e\neq e'}\|\Delta_e
\Delta_{e'}\R\|_2^2&\leq& 2\sup
_{e'}\bigl\{ A\bigl(e'\bigr)B
\bigl(e'\bigr)\bigr\}\sum_eA(e)B(e)
\EE\biggl[\frac{i_r^4(e)}{r^2(e)}\biggr]
\\
&\leq& 2 \sup_{e}\bigl\{A^2(e)B^2(e)
\bigr\}\sum_eA(e)B(e)C(e)\EE\bigl[r^2(e)i_r^4(e)
\bigr],
\end{eqnarray*}
where $C(e)$ is given by (\ref{eq:decoupling}):
\[
C(e)= \frac{\EE[{1}/({r^2(e)})\max\{1,{1}/({r^4(e)})\}]}{\EE[\min\{
r^2(e),{1}/({r^2(e)})\}]}.
\]
Thus, using Lemma~\ref{lem:thetae4Deltae},
where
\[
D(\mu)=2\sup_{e}\bigl\{A^2(e)B^2(e)
\bigr\}\sup_e\frac{A(e)B(e)C(e)}{\alpha_-(e)}.
\]


Finally, elementary calculus shows that
\[
D(\mu)\leq C\sup_e\bigl(\EE\bigl[r^8(e)\bigr]+
\EE\bigl[r^{-8}(e)\bigr]\bigr)^6\sup_e
\frac{\EE
[(r(e)-r'(e))^2 ({1}/({r^6(e)})+r^6(e) ) ]}{\EE
[(r(e)-r'(e))_+^2\min \{1,{1}/({r^4(e)}) \} ]},
\]
for some universal constant $C$.
\end{pf}

\begin{rem}
$K(\mu)$ is finite if the resistances are elliptic, or alternatively if
\[
\sup_e\bigl\{\EE\bigl(r^8(e)\bigr)+\EE
\bigl(r^{-8}(e)\bigr)\bigr\}<\infty
\]
and
\[
\inf_e\EE\bigl|r(e)-\EE\bigl(r(e)\bigr)\bigr|>0.
\]
Indeed, in this case, Cauchy--Schwarz inequality shows that
\[
\inf_e\EE \biggl[\bigl(r(e)-r'(e)
\bigr)_+^2\min \biggl\{1,\frac{1}{r^4(e)} \biggr\} \biggr]<0.
\]
\end{rem}

An easy consequence is that point-to-point effective resistances are
uniformly stable to noise.

\begin{coro}
\label{coro:noise}
On any random network with independent uniformly elliptic resistances,
or more generally with $K(\mu)$ finite, point-to-point effective
resistances are uniformly stable to noise:
\[
\inf_{u,v\in\V}\Corr\bigl(\R^{u,v}(r),\R^{u,v}
\bigl(r^\eps\bigr)\bigr)\mathop{\longrightarrow}_{\eps\rightarrow0}1.
\]
\end{coro}

\begin{pf}
Let us fix $u$ and $v$, and let $f$ denote the function $r\mapsto\R
^{u,v}(r)$. From Proposition~\ref{prop:noise},
\begin{eqnarray*}
\Cov\bigl(f(r),f\bigl(r^\eps\bigr)\bigr)&=&\sum
_{S\neq\varnothing}(1-\eps)^{|S|}\|f_S
\|^2
\\
&=&\Var(f)-\sum_{S\neq\varnothing}\bigl(1-(1-
\eps)^{|S|}\bigr)\|f_S\|^2
\\
&\geq&\Var(f)-\log\frac{1}{1-\eps}\sum_{S\neq\varnothing}|S|
\|f_S\|^2
\\
&\geq&\Var(f) \biggl(1-K(\mu)\log\frac{1}{1-\eps} \biggr).
\end{eqnarray*}
This shows that $\Corr(\R^{u,v}(r),\R^{u,v}(r^\eps))$ tends uniformly
(in $u$ and $v$) to 1 when $\eps$ tends to zero.
\end{pf}
A trivial consequence of Theorem~\ref{theo:EStight} is that the
Efron--Stein inequality, Corollary~\ref{coro:ESI}, is always tight for
point-to-point effective resistances.

\begin{coro}
\label{coro:EStight}
For any graph $G$ and any pair of vertices $(u,v)$,
\[
\Var\bigl(\R^{u,v}\bigr)\leq\sum_{e\in\E_{1/2}}\bigl\|
\Delta_e\R^{u,v}\bigr\|_2^2\leq K(\mu)\Var
\bigl(\R^{u,v}\bigr).
\]
\end{coro}

To finish this section, we shall recall the setting of the effective
resistance through a box in $\ZZ^d$, which was described in the
\hyperref[sec1]{Introduction}. Let $\B_n=\{0,\ldots,n\}^d$ be equipped with random
resistances $r$ on its set of edges $\E_n$, and let us define $A_n=\{
x\in\B_n\mbox{ s.t. }x_1=0\}$ and $Z_n=\{x\in\B_n\mbox{ s.t. }x_1=n\}$.
Then consider the graph with vertex set $B_n$, where all the
vertices of $A_n$ are identified to a single vertex on one side, and
all the vertices of $Z_n$ are identified on the other side. The \textit{effective resistance} through the box $\B_n$ is defined as
\[
\R_n=\R^{A_n,Z_n}_r.
\]
In \cite{Wehr97}, it is shown that under some hypotheses on the
distribution of the resistances,
%
\begin{equation}
\label{eq:Wehr}
\Var(\C_n)\geq Cn^{d-4},
\end{equation}
where $\C_n=1/\R_n$, and $C$ is some positive constant, depending on
$d$ and the common distribution of the conductances. Let us show how
one can recover (\ref{eq:Wehr}) in our setting. One could work directly
on $\C_n$, but will we rather show that
%
\begin{equation}
\label{eq:Wehrbis}
\Var(\R_n)\geq Cn^{4-3d},
\end{equation}
and then translate this bound on $\Var(\C_n)$. Let us suppose that
%
\begin{equation}
\label{eq:condmomentmean}
\sup_n\sup_{e\in\E_n}\bigl\{\EE
\bigl[r(e)\bigr]+\EE \bigl[r^{-1}(e)\bigr]\bigr\}<\infty.
\end{equation}
Then one may see that $\EE[\R_n]=\Theta(n^{2-d})$, and $\EE[\C_n]=\Theta
(n^{d-2})$ (cf. the argument in Section~5 of \cite{BenjaminiRossignol07}).


Let us call $i_{r}$ the unit current in the definition of $\R
^{A_n,Z_n}_r$, and order the edges in $E_n$ in some arbitrary fixed
way. We have the following martingale representation:
\begin{eqnarray*}
\R_n-\EE[\R_n]&=&\sum_{e\in\E_n}
\EE\bigl[\R_n(r)|(r_{e'})_{e'\geq e}\bigr]-\EE\bigl[\R
_n(r)|(r_{e'})_{e'> e}\bigr]
\\
&=&\sum_{e\in\E_n}\EE\bigl[\Delta_e
\R_n(r)|(r_{e'})_{e'\geq e}\bigr].
\end{eqnarray*}
Thus,
\begin{eqnarray*}
\Var(\R_n)&= &\sum_{e\in\E_n}\EE \bigl[\EE
\bigl[\Delta_e\R _n(r)|(r_{e'})_{e'\geq e}
\bigr]^2 \bigr]
\\
&=&\frac{1}{2} \sum_{e\in\E_n}\EE \bigl[\EE\bigl[
\R_n(r)-\R_n\bigl(r^{e\leftarrow
r'}\bigr)|(r_{e'})_{e'\geq e},r'(e)
\bigr]^2 \bigr]
\\
&=&\frac{1}{2} \sum_{e\in\E_n}\EE \bigl[\EE\bigl[
\R_n(r)-\R_n\bigl(r^{e\leftarrow
r'}\bigr)_+|(r_{e'})_{e'\geq e},r'(e)
\bigr]^2 \bigr]
\\
&&{}+\frac{1}{2} \sum_{e\in\E_n}\EE \bigl[\EE\bigl[
\R_n(r)-\R_n\bigl(r^{e\leftarrow
r'}\bigr)_-|(r_{e'})_{e'\geq e},r'(e)
\bigr]^2 \bigr],
\end{eqnarray*}
since the sign of $(\R_n(r)-\R_n(r^{e\leftarrow r'}))$ is the sign of
$r(e)-r'(e)$. Using inequality~(\ref{eq:majoDeltaeR2sided}),
\begin{eqnarray*}
\Var(\R_n)&\geq&\frac{1}{2} \sum
_{e\in\E_n}\EE \bigl[\bigl(r(e)-r'(e)
\bigr)_+^2\EE \bigl[i_r^2(e)|r(e),r'(e)
\bigr]^2 \bigr]
\\
&&{}+\frac{1}{2} \sum_{e\in\E_n}\EE \bigl[
\bigl(r(e)-r'(e)\bigr)_-^2\EE \bigl[i_{r^{e\leftarrow r'}}^2(e)|r(e),r'(e)
\bigr]^2 \bigr]
\\
& = &\sum_{e\in\E_n}\EE \bigl[\bigl(r(e)-r'(e)
\bigr)^2\EE\bigl[i_r^2(e)|r(e)
\bigr]^2 \bigr]
\\
& \geq &\sum_{e\in\E_n}\EE \biggl[\bigl(r(e)-r'(e)
\bigr)^2\min\biggl\{1,\frac
{1}{r^4(e)}\biggr\}\EE
\bigl[i_{r^{e\leftarrow1}}^2(e)|r(e)\bigr]^2 \biggr]
\\
& = &\sum_{e\in\E_n}\EE \biggl[\bigl(r(e)-r'(e)
\bigr)^2\min\biggl\{1,\frac{1}{r^4(e)}\biggr\} \biggr]\EE
\bigl[i_{r^{e\leftarrow1}}^2(e)\bigr]^2
\\
& \geq&\sum_{e\in\E_n}\frac{\EE [(r(e)-r'(e))^2\min\{1,{1}/({r^4(e))}\} ]}{\EE[\max\{r(e),{1}/({r(e)})\}]^2}\EE
\bigl[r(e)i_r^2(e)\bigr]^2,
\end{eqnarray*}
where the inequalities follows from Lemma~\ref{lem:regularite2} as in
the proof of~(\ref{eq:decoupling}). Define
\[
\tilde\alpha:=\inf_n\inf_{e\in E_n}
\frac{\EE [(r(e)-r'(e))^2\min\{
1,{1}/({r^4(e)})\} ]}{\EE[\max\{r(e),{1}/({r(e)})\}]^2}.
\]
Then, using Jensen's inequality,
\begin{eqnarray*}
\Var(\R_n)&\geq&\tilde\alpha\sum_{e\in\E_n}
\EE\bigl[r(e)i_r^2(e)\bigr]^2
\\
& \geq&\tilde\alpha\# (\E_n) \biggl(\frac{1}{\# (\E_n)}\sum
_{e\in\E
_n}\EE\bigl[r(e)i_{r}^2(e)\bigr]
\biggr)^2
\\
& = &\tilde\alpha\frac{1}{\# (\E_n)}\EE[\R_n]^2
\\
&\geq&C\tilde\alpha n^{4-3d}.
\end{eqnarray*}
This shows inequality (\ref{eq:Wehrbis}) under the condition $\tilde
\alpha>0$. Now, suppose that for some $p>2$,
%
\begin{equation}
\label{eq:condmomentWehr} \sup_n\sup_{e\in\E_n}\EE
\bigl[r(e)^p\bigr]<\infty,
\end{equation}
and let us show how one may recover (\ref{eq:Wehr}) from (\ref
{eq:Wehrbis}). Notice first that bounding $\R_n$ by the energy of the
deterministic flow which splits uniformly on $A_n$ and goes straight to
$Z_n$, there is a constant $c_d$ such that
\[
\R_n\leq c_d \frac{\sum_{e\in\E_n}r(e)}{n^{2d-2}}=:S_n.
\]
Using Rosenthal's inequality (see \cite{Rosenthal70}, Theorem~3), there
is a finite positive constant $k(p)$ such that for any $c>1$ and $p>2$,
\begin{eqnarray*}
\PP\bigl(\R_n\geq c\EE[S_n]\bigr)&\leq&\PP
\bigl(S_n\geq c\EE[S_n]\bigr)
\\
&\leq& \frac{\EE[|S_n-\EE(S_n)|^p]}{(c-1)^p\EE[S_n]^p}
\\
&= & \frac{\EE[|\sum_{e\in\E_n}(r(e)-\EE[r(e)])|^p]}{(c-1)^p\EE[\sum_{e\in\E_n}r(e)]^p}
\\
&\leq& \frac{k(p)\max\{\sum_{e\in\E_n}\EE[r(e)^p], (\sum_{e\in\E
_n}\EE[r(e)^2] )^{1/2}\}}{(c-1)^p\EE[\sum_{e\in\E_n}r(e)]^p}
\\
&\leq& \frac{k'(p)n^{d(1-p)}}{(c-1)^p}.
\end{eqnarray*}
Now,
\begin{eqnarray*}
\Var(\C_n)&\geq&\EE \biggl[ \biggl(\frac{1}{\R_n}-
\frac{1}{\EE[\R_n]} \biggr)^2 \biggr]
\\
&=&\EE \biggl[\frac{(\R_n-\EE[\R_n])^2}{\R_n^2\EE[\R_n]^2} \biggr]
\\
&\geq&\frac{\Var(\R_n)}{c^2\EE[\R_n]^2\EE[S_n]^2}-\EE \biggl[\frac{(\R
_n-\EE[\R_n])^2}{\R_n^2\EE[\R_n]^2}\II_{\R_n >c\EE[S_n]}
\biggr]
\\
&\geq&C\tilde\alpha\frac{n^{4-3d}}{c^2n^{8-4d}}-\EE \biggl[\frac{(\R
_n-\EE[\R_n])^2}{\R_n^2\EE[\R_n]^2}
\II_{\R_n >c\EE[S_n]} \biggr]
\\
&\geq&C\tilde\alpha\frac{n^{d-4}}{c^2}-\frac{\PP(\R_n\geq c\EE
[S_n])}{\EE[\R_n]^2}
\\
&\geq&C\tilde\alpha\frac{n^{d-4}}{c^2}-C'n^{2d-4}
\frac
{k'(p)n^dn^{d(1-p)}}{(c-1)^p}
\\
&\geq&C'n^{d-4},
\end{eqnarray*}
for $n$ large enough, since $p>2$. This shows inequality (\ref
{eq:Wehr}) when $\tilde\alpha$ is finite and (\ref{eq:condmomentmean})
and (\ref{eq:condmomentWehr}) hold. To state a simple moment condition,
for i.i.d. resistances with positive variance, one gets (\ref{eq:Wehr})
under the condition that $\EE[r^p(e)]+\EE[c(e)]<\infty$ for some $p>2$.
Notice that the moments of order $p>2$ are used only to go from the
bound on $\Var(\R_n)$ to a bound on $\Var(\C_n)$. Alternatively, one
could work directly on $\C_n$, thanks to the formula~(\ref
{eq:conductancemin}). One would need (this is not difficult) to
establish the analogs of Lemma~\ref{lem:regularite1} and \ref
{lem:regularite2} for $v_r$, the minimizer in (\ref
{eq:conductancemin}), instead of $i_r$. Then the line of proof which
gave (\ref{eq:Wehrbis}) would lead to (\ref{eq:Wehr}). To state a
simple moment condition, for i.i.d. conductances with positive
variance, one would obtain~(\ref{eq:Wehr}) under the condition that $\EE
[r(e)]+\EE[c(e)^2]<\infty$. In any case, our conditions seem to be much
weaker than the conditions in \cite{Wehr97}. For instance, no power-law
distribution satisfies Wehr's assumption, and he requires absolute
continuity w.r.t. the Lebesgue measure. We record the result of our
calculations in the following lemma, for the neat case where the
resistances are i.i.d.

\begin{lemm}
Suppose that the resistances are i.i.d., not constant and that $r(e)$
and $c(e)$ are integrable. Then
\[
\Var(\R_n)\geq Cn^{4-3d}.
\]
If furthermore $r(e)$ has a finite moment of order $p>2$, then
\[
\Var(\C_n)\geq C'n^{d-4}.
\]
The positive constants $C$ and $C'$ depend only on the common
distribution of the resistances.
\end{lemm}

Finally, let us emphasize the fact that we are unfortunately unable to
show that
%
\begin{equation}
\label{eq:varoptcube} \EE \biggl[\sum_{e\in\E_n}r^2(e)i_{r,n}^4(e)
\biggr]\leq Cn^{4-3d}.
\end{equation}
Otherwise, we would obtain from Corollary~\ref{coro:EStight} the
correct order for the variance of $\R_n$. Notice that (\ref
{eq:varoptcube}) does not necessarily hold when the resistances are not
supposed to have identical distribution, even in the elliptic case. To
understand why, let $v_r$ be the discrete harmonic function on $\B
_n\setminus(A_n\cup Z_n)$ with value $1$ on $A_n$ and $0$ on $Z_n$.
Discrete harmonic at $x$ means that
\[
d^*(c.dv_r) (x)=0.
\]
Then (\ref{eq:varoptcube}) is equivalent to
%
\begin{equation}
\label{eq:varoptcube2} \EE \biggl[\frac{1}{\# \E_n}\sum_{e\in\E_n}
\bigl(ndv_r(e)\bigr)^4 \biggr]\leq C'.
\end{equation}
When $n$ goes to infinity, one may compare $v_r$ to a continuous
analog. Let $(c_x)_{x\in\RR^d}$ be a deterministic elliptic collection
of conductance matrices. Let $\tilde v$ be the function on $[0,1]^d$
with value 1 on $\{x\mbox{ s.t. }x_1=0\}$, 0 on $\{x\mbox{ s.t. }x_1=1\}
$ and satisfying on $(0,1)^d$:
\[
\operatorname{div}(c.\nabla\tilde v)=0.
\]
The continuous analog of (\ref{eq:varoptcube2}) is the fact that
$\nabla\tilde v$ belongs to $L^4([0,1]^d)$. However, it is well known
that this may be false if the ellipticity constant $\Lambda$ is not
close enough to 1. On $\RR^2$, a counterexample is given in \cite
{DelmotteDeuschel05}; see the discussion after Proposition~1.1 therein.

\subsection{Further results for elliptic networks with homogeneous currents}
\label{subsec:homogeneous}

\textit{From now on}, \textit{the networks will be elliptic}: $r$ \textit{belongs\vspace*{1pt}
to}
$[1,\Lambda]^{\E_{1/2}}$ \textit{for some} $\Lambda\geq1$.
Recall that all
sets $\ell^2_-(\E,r)$ are the same for $r$ in $[1,\Lambda]^{\E_{1/2}}$,
since the norms with seights $r$ are all equivalent, the common set is
denoted by $\ell^2_-(\E)$; cf.~(\ref{eq:defellE}) and we shall refer to
the common norm topology of the sets $\ell^2_-(\E,r)$ as the
\textit{strong topology}. We shall consider the \textit{graph distance}, denoted
by $d$, on $\V$. Then, if $e$ and $e'$ are two edges in~$\E$, let
$d(e,e')$ be the maximal distance between two endpoints, one of which
is in $e$ and the other in $e'$. For any\vspace*{1pt} edge $e$ and any collection of
resistances $r$ in $[1,\Lambda]^{\E_{1/2}}$, the flow $i^e_r$ belongs
to $\ell^2_-(\E,r)$. Thus,
%
\begin{equation}
\label{eq:energylocalized}
\sum_{e'\dvtx  d(e',e)\geq L}r\bigl(e'
\bigr) \bigl(i^e_r\bigl(e'\bigr)
\bigr)^2\mathop{\longrightarrow}_{L\rightarrow\infty}0.
\end{equation}
Below, we shall be interested in graphs where the above convergence
holds \textit{uniformly in $e$ and $r$}. We shall say that such a graph
has \textit{homogeneous currents}.

\begin{defi}
Let $\G=(\V,\E)$ be a countable, oriented, symmetric and connected
graph. Let $\Lambda\geq1$ be a real number, and define
%
\begin{equation}
\label{eq:alphaG}
\hspace*{6pt}\alpha(\G,L,\Lambda)=\sup \biggl\lbrace\sum
_{e'\dvtx  d(e',e)\geq
L}r\bigl(e'\bigr) \bigl(i^e_r
\bigl(e'\bigr)\bigr)^2\mbox{ s.t. }e\in\E, r\in[1,
\Lambda]^{\E
_{1/2}} \biggr\rbrace.
\end{equation}
The graph $\G$ is said to have $\Lambda$-\textit{homogeneous currents} if
\[
\alpha(\G,L,\Lambda)\mathop{\longrightarrow}_{L\rightarrow+\infty}0.
\]
\end{defi}

It is natural to expect that for every $\Lambda\geq\Lambda'$ strictly
larger than 1, $\G$ has $\Lambda$-homogeneous currents if it has
$\Lambda'$-homogeneous currents [the other direction being trivial
since $\alpha(\G,L,\Lambda)$ is monotone in $\Lambda$]. However, we
could not prove this.

The first fundamental observation, due to Mika\"el de la Salle, is that
for a fixed edge $e$, the convergence in (\ref{eq:energylocalized})
always hold uniformly in $r$, thanks to a compactness argument.

\begin{prop}
\label{prop:Mikael}
Let $u$ and $v$ be two vertices of $\G$, and suppose that $(\G
_L)_{L\geq0}$ is a sequence of finite connected graphs that exhausts
$\G$ and such that~$\G_0$ contains $u$ and $v$. Then
\[
\sup \biggl\lbrace\sum_{e'\in\G_L^c}r\bigl(e'
\bigr) \bigl(i^{u,v}_r\bigl(e'\bigr)
\bigr)^2\mbox{ s.t. }r\in[1,\Lambda]^{\E_{1/2}} \biggr\rbrace
\mathop{\longrightarrow}_{L\rightarrow +\infty}0.
\]
\end{prop}

\begin{pf*}{Proof (Due to Mika\"el de la Salle)}
Let us fix $u$ and $v$ two vertices of~$\G$. We equip $[1,\Lambda]^{\E
_{1/2}}$ with the product topology. For a fixed $\theta\in\ell_-^2(\E)$
and $\eps>0$, let $F$ be a finite subset of edges such that $\sum_{e\in
F^c}\theta^2(e)<\eps$. Then
\[
\biggl|\sum_{e}r(e)\theta^2(e)-\sum
_{e}r'(e)\theta^2(e)\biggr|\leq\sum
_{e\in
F}\bigl|r(e)-r'(e)\bigr|\theta^2(e)+
\Lambda\eps,
\]
and thus the function $r\mapsto\sum_{e}r(e)\theta^2(e)$ is continuous
for the product topology. Then $r\mapsto\R^{u,v}(r)$ is an infimum of
continuous functions on $[1,\Lambda]^{\E_{1/2}}$, and thus, it is upper
semi-continuous on $[1,\Lambda]^{\E_{1/2}}$. Now, define
\[
c:=\lim_{L\rightarrow+\infty} \sup \biggl\lbrace\sum
_{e'\in\G
_L^c}r\bigl(e'\bigr) \bigl(i^{u,v}_r
\bigl(e'\bigr)\bigr)^2\mbox{ s.t. }r\in[1,
\Lambda]^E \biggr\rbrace,
\]
which exists by the monotonicity in $L$ of the\vspace*{1pt} right-hand side.

One may find a sequence $(r_L)_{L\geq1}$ in $[1,\Lambda]^{\E_{1/2}}$
such that
\[
c=\lim_{L\rightarrow+\infty} \sum_{e'\in\G
_L^c}r_L
\bigl(e'\bigr) \bigl(i^{u,v}_{r_L}
\bigl(e'\bigr)\bigr)^2.
\]
By Lemma~\ref{lem:boundcurrent}, the sequence $(i^{u,v}_{r_L})_{L\geq
1}$ lies in $[-1,1]^{\E}$ [and even in $\ell^2_-(\E)$] which is compact
for the product topology by Cantor's diagonal argument. Also, the
sequence $(r_L)_{L\geq1}$ lies in the compact set $[1,\Lambda]^{\E
_{1/2}}$. Thus, one may find an increasing sequence of integers
$(L_k)_{k\geq1}$ such that for any $e\in\E$, $(r_{L_k})_{k\geq1}(e)$
converges to some value $r(e)$ in $[1,\Lambda]$ and
$(i^{u,v}_{r_{L_k}}(e))_{k\geq1}$ converges to some $\theta(e)\in
[-1,1]$. Then $\theta$ is a unit flow from $u$ to $v$. From the
upper-semi-continuity of $r\mapsto\R^{u,v}(r)$,
\begin{eqnarray*}
\R^{u,v}(r)&\geq&\limsup_{k\rightarrow+\infty} \R^{u,v}(r_{L_k})
\\
&=&\limsup_{k\rightarrow+\infty}\sum_{e'}r_{L_k}
\bigl(e'\bigr) \bigl(i^{u,v}_{r_{L_k}}
\bigl(e'\bigr)\bigr)^2
\\
&=&\limsup_{k\rightarrow+\infty} \biggl(\sum_{e'\in\G
_{L_k}}r_{L_k}
\bigl(e'\bigr) \bigl(i^{u,v}_{r_{L_k}}
\bigl(e'\bigr)\bigr)^2+\sum
_{e'\in\G
_{L_k}^c}r_{L_k}\bigl(e'\bigr)
\bigl(i^{u,v}_{r_{L_k}}\bigl(e'\bigr)
\bigr)^2 \biggr)
\\
&= &\limsup_{k\rightarrow+\infty}\sum_{e'\in\G
_{L_k}}r_{L_k}
\bigl(e'\bigr) \bigl(i^{u,v}_{r_{L_k}}
\bigl(e'\bigr)\bigr)^2+c
\\
&\geq&\sup_{L'}\limsup_{k\rightarrow+\infty}\sum
_{e'\in\G
_{L'}}r_{L_k}\bigl(e'\bigr)
\bigl(i^{u,v}_{r_{L_k}}\bigl(e'\bigr)
\bigr)^2+c
\\
&= &\sup_{L'}\sum_{e'\in\G_{L'}}r
\bigl(e'\bigr) \bigl(\theta\bigl(e'\bigr)
\bigr)^2+c
\\
&=&\sum_{e'}r\bigl(e'\bigr) \bigl(
\theta\bigl(e'\bigr)\bigr)^2+c
\\
&\geq& \R^{u,v}(r)+c.
\end{eqnarray*}
This shows that $c=0$ and proves the proposition.
\end{pf*}
Notice that in Proposition~\ref{prop:Mikael}, the ellipticity
hypothesis is used in a crucial way. This proposition will allow us to
find our first graphs with homogeneous currents: the
\textit{quasi-transitive} graphs. There is no univeral definition, but in this
article we shall say that a graph $\G=(\V,\E)$ with automorphism group
$\operatorname{\mathsf{Aut}}(G)$ is quasi-transitive if its set of edges $\E$ is
composed of a finite number of distinct orbits under the natural action
of $\operatorname{\mathsf{Aut}}(G)$ on  $\E$.

\begin{coro}
\label{coro:transitivehomogeneous}
Let $\G=(\V,\E)$ be a countable, oriented, symmetric and connected
graph. Suppose that $\G$ is quasi-transitive. Then, for any $\Lambda
\geq1$, $\G$ has \textup{$\Lambda$-homogeneous currents}.
\end{coro}

\begin{pf}
The quasi-transitivity hypothesis implies that there exists a finite
set of edges $e_1,\ldots,e_r$ such that
\[
\alpha(\G,L,\Lambda)=\max_{i=1,\ldots,r} \sup \biggl\lbrace\sum
_{e'\dvtx d(e',e_i)\geq L}r\bigl(e'\bigr)
\bigl(i^{e_i}_r\bigl(e'\bigr)
\bigr)^2\mbox{ s.t. }r\in[1,\Lambda]^{\E
_{1/2}} \biggr\rbrace.
\]
But, from Proposition~\ref{prop:Mikael}, for any $i$,
\[
\sup \biggl\lbrace\sum_{e'\dvtx  d(e',e_i)\geq
L}r\bigl(e'
\bigr) \bigl(i^{e_i}_r\bigl(e'\bigr)
\bigr)^2\mbox{ s.t. }r\in[1,\Lambda]^{\E_{1/2}} \biggr\rbrace
\mathop{\longrightarrow}_{L\rightarrow+\infty}0,
\]
taking $\G_L$ to be the graph whose edges are all the edges of $\E$ at
distance at most~$L$ from $e_i$ and whose vertices are the endpoints of
those edges. Thus, $\alpha(\G,L,\Lambda)$ goes to zero as $L$ goes to
infinity, and $\G$ has $\Lambda$-homogeneous currents.
\end{pf}

A consequence of Corollary~\ref{coro:transitivehomogeneous} is that $\ZZ
^d$ has homogeneous currents. This can also be seen in a more robust
way using the powerful tool of \textit{elliptic Harnack inequality}.
Indeed, let $B_L(e)$ be the vertices at distance at most $L$ from $e$.
Then $\sum_{d(e',e)\geq L}r(e')(i^e_r(e'))^2$ is upper-bounded by the
oscillation on $B_L(e)^c$ of the votage induced\vspace*{1pt} by the flow $i^e_r$.
Since this voltage is a bounded function, harmonic on $\ZZ^d\setminus
e$, one may then show, using the Harnack inequality of \cite{DelmotteElliptic} as in~\cite{MoserElliptic}, Section~6, that $\alpha
(\ZZ^d,L,\Lambda)$ decays at least as quickly as a negative power of
$L$. This argument can be carried out on any graph satisfying the
conditions of~\cite{DelmotteElliptic} plus the additional condition
that the annuli between $B_L(e)$ and $B_{4L}(e)$ are connected and may
be covered by a bounded number of balls of radius $L$. For instance,
this shows that any graph roughly isometric to $\ZZ^d$ has homogeneous currents.

\begin{figure}

\includegraphics{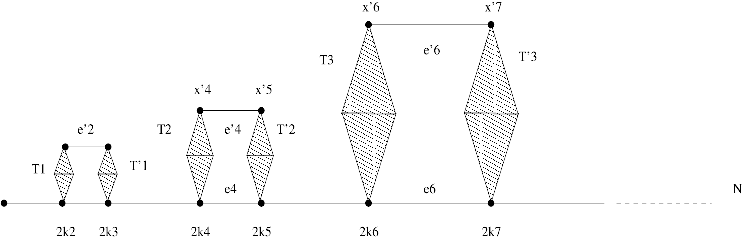}

\caption{An example of a graph without homogeneous currents.}
\label{fig:nonhomogeneous}
\end{figure}

Now, let us give an (artificial) example of a graph which does not have
homogeneous currents. A perfect binary tree of depth $k$ is a rooted
binary tree where every leaf is at depth $k$ and all other vertices
have two children. For any $k\in\NN^*$, let $T_k$ be two copies of a
perfect binary trees of depth $k$ glued at the leafs. The result has
two roots. Now, to construct our graph (see Figure~\ref{fig:nonhomogeneous}), we start from $\NN^*$, with the usual notion of
graph on it, and for any $k\in\NN^*$, we do the following construction.
We add $T_k$ by glueing one of its roots on the vertex $x_{2k}:=2k$ of
$\NN$ and call the remaining root $x'_{2k}$. Then we add a copy $T'_k$
of $T_k$ by gluing one of its root on the vertex $x_{2k+1}:=2k+1$ of
$\NN$ and call the remaining root $x'_{2k+1}$. Then we join $x'_{2k}$
and $x'_{2k+1}$ by an edge.\vspace*{1pt} We denote by $\G$ the resulting graph, let
$e_{2k}:=(x_{2k},x_{2k+1})$, $e'_{2k}:=(x'_{2k},x'_{2k+1})$ and equip
this graph with unit resistances. Now,\vspace*{-1pt} the resistance between the root
and the leafs of a complete binary tree of any depth is at most~$1/2$.
Thus,\vspace*{-1pt} $\R^{x_{2k},x'_{2k}}_{T_k}$, the resistance between $x_{2k}$ and
$x'_{2k}$ in the graph $T_k$ is at most~$1$. Similarly, $\R
^{x_{2k+1},x'_{2k+1}}_{T'_k}$ is at most~1. Since resistances in series
add, one sees that
\[
i^{e_{2k}}_r\bigl(e'_{2k}\bigr)=
\frac{r(e_{2k})}{\R
^{x_{2k},x'_{2k}}_{T_k}+r(e'_{2k})+\R^{x_{2k+1},x'_{2k+1}}_{T'_k}}\geq \frac{1}{3}.
\]
On the other hand, the graph-theoretical distance between $e_{2k}$ and
$e'_{2k}$ is $2k+1$. Thus, for any $\Lambda\geq1$, and $L\geq1$,
\[
\alpha(\G,L,\Lambda)\geq\sum_{e'\dvtx  d(e',e_{2L})\geq
L}r(e)
\bigl(i_r^{e_{2L}}\bigl(e'\bigr)
\bigr)^2\geq\bigl(i_r^{e_{2L}}
\bigl(e'_{2L}\bigr)\bigr)^2\geq
\frac
{1}{9}.
\]
Whence $\G$ does not have homogeneous currents.

Finally, the reasoning in the proof of Proposition~\ref{prop:Mikael}
allows also to obtain some regularity of the functions $r\mapsto
i^{u,v}_r$ and $r\mapsto\R^{u,v}(r)$ when $[1,\Lambda]^{\E_{1/2}}$ is
equipped with the product topology (it will not be used in the sequel).

\begin{prop}
\label{prop:continuity}
Let $u$ and $v$ be two vertices of $\G$. If $[1,\Lambda]^{\E_{1/2}}$ is
equipped with the product topology and $\ell_-^2(\E)$ with the strong
topology, then the following maps are continuous:
\[
\cases{
\ds [1,\Lambda]^{\E_{1/2}}\rightarrow
\ell_-^2(\E),
\vspace*{2pt}\cr
\ds r \mapsto  i_r^{u,v},}
\quad\mbox{and}\quad
\cases{
\ds [1,\Lambda]^{\E_{1/2}} \rightarrow \RR,
\vspace*{2pt}\cr
r  \mapsto   \R^{u,v}(r).}
\]
\end{prop}

\begin{pf}
Take any $r$ in $[1,\Lambda]^E$ and any sequence $(r_L)$ converging to
$r$. As in the proof of Proposition~\ref{prop:Mikael}, one may extract
a sequence $(r_{L_k})_{k\geq1}$ such that $(i^{u,v}_{r_{L_k}})_{k\geq
1}$ converges pointwise to some $\theta$ which is a unit flow from $u$
to $v$. Then
\begin{eqnarray*}
\R^{u,v}(r) &\geq & \limsup_{k\rightarrow+\infty} \R^{u,v}(r_{L_k})
\geq \liminf_{k\rightarrow+\infty} \R^{u,v}(r_{L_k})
\\
&\geq&\liminf_{k\rightarrow+\infty}\sum_{e'\in\G
_{L_k}}r_{L_k}
\bigl(e'\bigr) \bigl(i^{u,v}_{r_{L_k}}
\bigl(e'\bigr)\bigr)^2
\\
&\geq&\sum_{e'}r
\bigl(e'\bigr) \bigl(\theta\bigl(e'\bigr)
\bigr)^2,
\end{eqnarray*}
by Fatou's lemma. Since $i^{u,v}_{r}$ is the unique minimizer of $\theta
\mapsto\sum_{e'}r(e')\theta(e')$ over the unit flows from $u$ to $v$,
this shows at once that $\theta=i^{u,v}_{r}$ and that $\R
^{u,v}(r_{L_k})$ converges to $\R^{u,v}(r)$ as $k$ goes to infinity.
Since this is true for any subsequence $(r_{L_k})_k$ such that
$(i^{u,v}_{r_{L_k}})_{k}$ converges pointwise, we deduce that
$i^{u,v}_{r_L}$ converges pointwise to $i^{u,v}_{r}$ and that $\R
^{u,v}(r_{L})$ converges to $\R^{u,v}(r)$. Notably, $r\mapsto\R
^{u,v}(r)$ is continuous for the product topology.

Now, recall that $i_{r_L}$ is a current, thus there exists a function
$v$ such that $r(e)i_{r_L}(e)=dv(e)$ for any $e$. Notice also that
$d^*(i_{r_L}^{u,v}-i_r^{u,v})=0$. This implies that $i_{r_L}$ is
orthogonal in $\ell_-^2(\E,r)$ to $i_{r_L}^{u,v}-i_r^{u,v}$. Indeed,
suppose first that the network is finite. Then
\begin{eqnarray*}
\sum_er_L(e)i^{u,v}_{r_L}(e)
\bigl(i_{r_L}^{u,v}(e)-i_r^{u,v}(e)
\bigr)&=&\sum_edv(e) \bigl(i_{r_L}^{u,v}(e)-i_r^{u,v}(e)
\bigr)
\\
&=&\sum_{x\in\V}v(x)d^*\bigl(i_{r_L}^{u,v}-i_r^{u,v}
\bigr) (x)
\\
&=&0.
\end{eqnarray*}
This continues to hold when $\G$ is not finite since $i_{r_L}^{u,v}$
and $i_r^{u,v}$ are limits in $\ell_-^2(\E)$ (for the strong topology)
of wired currents on finite graphs. Thus,
\[
\bigl\|i^{u,v}_{r_L}-i^{u,v}_{r}
\bigr\|_{r_L}^2=\sum_{e'\in\E
_{1/2}}r_L
\bigl(e'\bigr) \bigl(i^{u,v}_{r}
\bigl(e'\bigr)\bigr)^2-\R^{u,v}(r_L).
\]
From the dominated convergence theorem, since $r_L$ converges pointwise
to $r$,
\[
\lim_{L\rightarrow\infty}\sum_{e'\in\E_{1/2}}r_L
\bigl(e'\bigr) \bigl(i^{u,v}_{r}
\bigr)^2=\R ^{u,v}(r).
\]
Thus,
\[
\lim_{L\rightarrow\infty}\bigl\|i^{u,v}_{r_L}-i^{u,v}_{r}
\bigr\|_{r_L}^2=0,
\]
which shows that $r\mapsto i^{u,v}_{r}$ is continuous.
\end{pf}

\subsubsection{Concentration on sets of small diameter}
A small variation on the proof of Theorem~\ref{theo:EStight} allows to
obtain the following result, which shows that on graphs with
homogeneous currents, the Walsh decomposition is concentrated on sets
of small diameter.

\begin{theo}
\label{theo:smalldiam}
For any graph $\G$, any $\Lambda\geq1$ and any $L\geq1$,
\[
\sum_{\operatorname{diam}(S)\geq L}\bigl\|\R^{u,v}_S
\bigr\|_2^2\leq C(\Lambda)\alpha(G,L,\Lambda )\sum
_{S\neq\varnothing}\bigl\|\R^{u,v}_S\bigr\|_2^2.
\]
\end{theo}

\begin{pf}
Let $L$ be a positive integer. From inequality (\ref{eq:MajoDeltaee'}),
one gets
\begin{eqnarray*}
&&\sum_{(e,e')\dvtx  d(e,e')\geq L}\bigl\|\Delta_e
\Delta_{e'}\R^{u,v}\bigr\|_2^2
\\
&&\qquad \leq C(\Lambda)\EE \biggl[\sum_{d(e,e')\geq
L}A(e)A
\bigl(e'\bigr) \bigl(i_r^4(e)+i_r^4
\bigl(e'\bigr)\bigr)\frac{i_r^e(e')^2}{r(e)^2}r(e)r\bigl(e'
\bigr) \biggr]
\\
&&\qquad =2C(\Lambda)\EE \biggl[\sum_{d(e,e')\geq L}A(e)A
\bigl(e'\bigr)i_r^4(e)
\frac
{i_r^e(e')^2}{r(e)^2}r(e)r\bigl(e'\bigr) \biggr]
\\
&&\qquad \leq 2C(\Lambda)\sup_{e'}A\bigl(e'\bigr)\EE
\biggl[\sum_e\frac
{1}{r(e)}A(e)i_r^4(e)
\mathop{\sum_{e' \mathrm{s.t.}}}_{d(e,e')\geq L}i_r^e
\bigl(e'\bigr)^2r\bigl(e'\bigr) \biggr]
\\
&&\qquad \leq 2C(\Lambda)\sup_{e'}A\bigl(e'\bigr)
\alpha(G,L,\Lambda)\EE\biggl[\sum_eA(e)r^2(e)i_r^4(e)
\biggr]
\\
&&\qquad \leq  C'(\Lambda)\alpha(G,L,\Lambda)\sum
_{S\neq\varnothing}\bigl\|\R ^{u,v}_S
\bigr\|_2^2.
\end{eqnarray*}
Then, Theorem~\ref{theo:EStight} leads to
\[
\sum_{d(e,e')\geq L}\bigl\|\Delta_e
\Delta_{e'}\R^{u,v}\bigr\|_2^2\leq C(
\Lambda )\alpha(G,L,\Lambda)\sum_{S\neq\varnothing}\bigl\|
\R^{u,v}_S\bigr\|_2^2.
\]
Finally, notice that
\begin{eqnarray*}
\sum_{d(e,e')\geq L}\bigl\|\Delta_e
\Delta_{e'}\R^{u,v}\bigr\|_2^2 &=& \sum
_{d(e,e')\geq L}\sum_{S\supset\{e,e'\}}\bigl\|
\R^{u,v}_S\bigr\|_2^2
\\
&=&\sum_{S}\mathop{\sum
_{e,e'\in S}}_{d(e,e')\geq L}\bigl\|\R^{u,v}_S\bigr\|
_2^2
\\
&\geq&\sum_{\operatorname{diam}(S)\geq L}\bigl\|\R^{u,v}_S
\bigr\|_2^2.
\end{eqnarray*}
\upqed\end{pf}
%
\subsubsection{The first level carries a significant weight}
A corollary of Theorem~\ref{theo:smalldiam} is that on graphs with
homogeneous currents, the first level of the Walsh decomposition
carries a significant part of the $L^2$-norm of the centered
point-to-point resistances.

\begin{coro}
\label{coro:firstlevel}
Suppose that $\G$ has $\Lambda$-homogeneous currents and degree bounded
by $\delta$. Then there is a constant $C(\Lambda,\delta,\G)$ such that
\[
\sum_{e}\bigl\|\R^{u,v}_{\{e\}}
\bigr\|_2^2\leq\Var\bigl(\R^{u,v}\bigr)=\sum
_{S\neq\varnothing}\bigl\|\R^{u,v}_S\bigr\|_2^2
\leq C(\Lambda,\delta,\G)\sum_{e}\bigl\|\R
^{u,v}_{\{e\}}\bigr\|_2^2.
\]
\end{coro}
\begin{pf}
Let us fix $u$ and $v$ two vertices of $G$, and let us drop the
superscript $u,v$ to lighten the notation. Since $G$ has $\Lambda
$-homogeneous currents, one may find $L$ large enough (depending on $G$
and $\Lambda$) so that
\[
\sum_{S\neq\varnothing}\|\R_S\|_2^2
\leq2\mathop{\sum_{\operatorname{diam}(S)\leq L}}_{S \neq\varnothing}\|
\R_S\|_2^2.
\]
For $L$ a positive integer, using Lemma~\ref{lemm:regmean} below,
\begin{eqnarray*}
&&\mathop{\sum_{\operatorname{diam}(S)\leq L}}_{S \neq\varnothing}\|
\R_S\|_2^2
\\
&&\qquad \leq  C(\Lambda,L)\sum_{\operatorname{diam}(S)\leq L}\sum
_{e\in S}\bigl\|\R_{\{e\}}\bigr\|_2^2
\\
&&\qquad  = C(\Lambda,L)\sum_{e\in\E}\bigl\|\R_{\{e\}}
\bigr\|_2^2\bigl|\bigl\{S \mbox{ s.t. } e\in S\mbox{ and }\operatorname{diam}(S)
\leq L\bigr\}\bigr|
\\
&&\qquad \leq  C'(\Lambda,L,\delta)\sum_{e\in\E}
\|\R_{\{e\}}\|_2^2.
\end{eqnarray*}
\upqed\end{pf}
In the proof above, we used the following lemma, which allows to
control $\|\R^{u,v}_{S}\|_p$ for small sets $S$. It will be used again
in the central limit approximation of Theorem~\ref{theo:TLCinfini}.

\begin{lemm}
\label{lemm:regmean}
For any $p\geq1$,
\[
C_1(\Lambda)m_p(e)\bigl\|\bigl(i^{u,v}_r(e)
\bigr)^2\bigr\|_1\leq\bigl\|\R^{u,v}_{\{e\}}
\bigr\|_p\leq C_2(\Lambda)m_p(e)\bigl\|
\bigl(i^{u,v}_r(e)\bigr)^2\bigr\|_1.
\]
Furthermore, for any $S$ such that $0<|S|\leq L$,
\[
\bigl\|\R^{u,v}_{S}\bigr\|_p\leq C(\Lambda)^L
\sum_{e\in S}\bigl\|\R^{u,v}_{\{e\}}\bigr\|
_p.
\]
\end{lemm}

\begin{pf}
Let $r$ and $r'$ be two independent random variables with the same
distribution $\PP:=\bigotimes_{e\in E_{1/2}}\mu_e$. For any function
$F$ in $L^2(\RR^{\E_{1/2}})$ and any $p\geq1$,
\begin{eqnarray*}
\|F_{\{e\}}\|_p^p&= & \EE \bigl[\bigl\llvert \EE
\bigl( \Delta_e F|r(e)\bigr)\bigr\rrvert^p \bigr]
\\
&=&\EE \bigl[\bigl\llvert \EE\bigl( F\bigl(r^{e\leftarrow r'}\bigr)-F(r)|r'(e)
\bigr)\bigr\rrvert ^p \bigr]
\\
&\leq&\EE \bigl[\bigl\llvert \EE\bigl( F\bigl(r^{e\leftarrow r'}
\bigr)-F(r)|r(e),r'(e)\bigr)\bigr\rrvert ^p \bigr]
\\
&=&2\EE \bigl[\bigl\llvert \EE \bigl( \bigl[F\bigl(r^{e\leftarrow r'}\bigr)-F(r)
\bigr]\II_{r'(e)>
r(e)}|r(e),r'(e) \bigr)\bigr\rrvert
^p \bigr].
\end{eqnarray*}
Recall from (\ref{eq:majoDeltaeR2sided}) that
\begin{eqnarray*}
0 &\leq & \bigl(r'(e)-r(e)\bigr)_+i_{r^{e\leftarrow r'}}^2(e)
\leq\bigl(\R\bigl(r^{e\leftarrow
r'}\bigr)-\R(r)\bigr)\II_{r'(e)> r(e)}\\
&\leq &
\bigl(r'(e)-r(e)\bigr)_+i_r^2(e).
\end{eqnarray*}
Lemma~\ref{lem:regularite2} allows then to decouple $(r'(e)-r(e))_+$
and $i_r^2(e)$ [or $(r'(e)-r(e))_+$ and $i_{r^{e\leftarrow r'}}^2(e)$]:
\begin{eqnarray*}
\|\R_{\{e\}}\|_p^p&\leq&2\EE\bigl[
\bigl(r'(e)-r(e)\bigr)_+^p\EE \bigl[i_r^2(e)|r(e),r'(e)
\bigr]^p\bigr]
\\
&\leq&C(\Lambda)^p\EE\bigl[\bigl(r'(e)-r(e)
\bigr)_+^p\EE\bigl[i_{r^{e\leftarrow
1}}^2(e)|r(e),r'(e)
\bigr]^p\bigr]
\\
& = &C(\Lambda)^p\EE\bigl[\bigl(r'(e)-r(e)
\bigr)_+^p\bigr]\EE\bigl[i_{r^{e\leftarrow
1}}^2(e)
\bigr]^p
\\
&\leq& C'(\Lambda)^pm_p(e)^p
\EE\bigl(i_r^2(e)\bigr)^p.
\end{eqnarray*}
On the other hand, for any function $F$ in $L^2(\RR^{\E_{1/2}})$ and
any $p\geq1$
\begin{eqnarray*}
&&\bigl\llvert \EE\bigl[F\bigl(r^{e\leftarrow r'}\bigr)-F(r)|r(e),r'(e)
\bigr]\bigr\rrvert^p
\\
&&\qquad \leq \bigl(\bigl\llvert \EE\bigl[F\bigl(r^{e\leftarrow r'}\bigr)-
\EE(F)|r(e),r'(e)\bigr]\bigr\rrvert +\bigl\llvert \EE\bigl[F(r)-
\EE(F)|r(e),r'(e)\bigr]\bigr\rrvert \bigr)^p
\\
&&\qquad \leq 2^{p-1} \bigl(\bigl\llvert \EE\bigl[F\bigl(r^{e\leftarrow r'}
\bigr)-\EE (F)|r(e),r'(e)\bigr]\bigr\rrvert ^p\\
&&\hspace*{21pt}\qquad\qquad{}+\bigl
\llvert \EE\bigl[F(r)-\EE(F)|r(e),r'(e)\bigr]\bigr\rrvert
^p \bigr),
\end{eqnarray*}
and thus:
\[
\EE \bigl[\bigl\llvert \EE\bigl[F\bigl(r^{e\leftarrow r'}\bigr)-F(r)|r(e),r'(e)
\bigr]\bigr\rrvert ^p \bigr]\leq2^p\EE \bigl[\bigl\llvert
\EE\bigl[F\bigl(r^{e\leftarrow r'}\bigr)-\EE (F)|r(e),r'(e)\bigr]\bigr
\rrvert ^p \bigr].
\]
Furthermore,
\[
F_{\{e\}}=\EE\bigl[F\bigl(r^{e\leftarrow r'}\bigr)-\EE(F)|r(e),r'(e)
\bigr]
\]
and
\begin{eqnarray*}
&&\EE \bigl[\bigl\llvert \EE\bigl[F\bigl(r^{e\leftarrow r'}\bigr)-F(r)|r(e),r'(e)
\bigr]\bigr\rrvert ^p \bigr]
\\
&&\qquad =2\EE \bigl[\bigl\llvert \EE\bigl[\bigl(F\bigl(r^{e\leftarrow r'}\bigr)-F(r)
\bigr)\II _{r(e)>r'(e)}|r(e),r'(e)\bigr]\bigr\rrvert
^p \bigr].
\end{eqnarray*}
Thus,
\[
\|F_{\{e\}}\|_p^p\geq\frac{1}{2^{p-1}}\EE
\bigl[\bigl\llvert \EE \bigl[\bigl(F\bigl(r^{e\leftarrow r'}\bigr)-F(r)\bigr)
\II_{r(e)>r'(e)}|r(e),r'(e)\bigr]\bigr\rrvert ^p
\bigr].
\]
%
Now, as for the upper bound, one uses (\ref{eq:majoDeltaeR2sided}) and
Lemma~\ref{lem:regularite2}:
\begin{eqnarray*}
\|\R_{\{e\}}\|_p^p&\geq&\frac{1}{2^{p-1}}\EE
\bigl[\EE\bigl[ \bigl(r'(e)-r(e)\bigr)_+i_{r^{e\leftarrow r'}}^2(e)|r(e),r'(e)
\bigr]^p \bigr]
\\
&\geq&C(\Lambda)^p\EE \bigl[\EE\bigl[ \bigl(r'(e)-r(e)
\bigr)_+i_{r^{e\leftarrow
1}}^2(e)|r(e),r'(e)
\bigr]^p \bigr]
\\
&=&C(\Lambda)^p\EE\bigl[\bigl(r'(e)-r(e)
\bigr)_+^p\bigr]\EE\bigl[i_{r^{e\leftarrow1}}^2(e)
\bigr]^p
\\
&\geq&m_p(e)^pC'(\Lambda)^p\EE
\bigl[i_r^2(e)\bigr]^p.
\end{eqnarray*}

This proves the first part of the lemma. Now, let $S$ be such that
$0<|S|\leq L$ and let $e\in S$. For a subset $S$ of $\E$, recall the
definition of $L_S$ in (\ref{eq:defLS}). Notice first that for any
function $F$ in $L^2(\RR^{\E_{1/2}})$,
\begin{eqnarray*}
F_S(r)&=&\EE \biggl[ \biggl(\prod_{e'\in S\setminus\{e\}}(1-L_{\{e\}
})
\biggr)\Delta_eF(r)\Big|r_S \biggr]
\\
&=&\EE \biggl[ \biggl(\sum_{I\subset S\setminus\{e\}
}(-1)^{|I|}L_I
\biggr)\Delta_eF(r)\Big|r_S \biggr]
\\
&=&\EE \biggl[\sum_{I\subset S\setminus\{e\}}(-1)^{|I|}L_I(
\Delta _eF) (r)\Big|r_S \biggr].
\end{eqnarray*}
Thus, using Jensen's inequality,
\[
\|F_{S}\|_p^p\leq2^{p|S\setminus\{e\}|}\EE \bigl[
\EE \bigl(\bigl|\Delta _eF(r)\bigr| |r_S \bigr)^p
\bigr].
\]
Since this is true for any $e\in S$, we have
%
\begin{equation}
\label{eq:FSFe}
\|F_{S}\|_p^p\leq\min
_{e\in S}2^{p|S\setminus\{e\}|}\EE \bigl[\EE \bigl(\bigl|
\Delta_eF(r)\bigr| |r_S \bigr)^p \bigr].
\end{equation}
Now,
\[
\EE \bigl(\bigl|\Delta_eF(r)\bigr| |r_S \bigr)^p
\leq\EE \bigl(\bigl|F(r)-F\bigl(r^{e\leftarrow r'}\bigr)\bigr| |r_S
\bigr)^p.
\]
Using (\ref{eq:majoDeltaeR2sided}) and Lemma~\ref{lem:regularite2},
\begin{eqnarray*}
\EE \bigl(\bigl|\R(r)-\R\bigl(r^{e\leftarrow r'}\bigr)\bigr| |r_S \bigr)&
\leq&C(\Lambda)\EE \bigl(\bigl|r(e)-r'(e)\bigr|i_{r^{e\leftarrow1}}^2(e)
|r_S \bigr)
\\
&= & C(\Lambda)\EE\bigl(\bigl|r(e)-r'(e)\bigr| |r_S\bigr)\EE
\bigl(i_{r^{e\leftarrow1}}^2(e) |r_S \bigr).
\end{eqnarray*}
Thus,
\begin{eqnarray*}
\EE \bigl[\EE \bigl(\bigl|\Delta_e\R(r)\bigr| |r_S
\bigr)^p \bigr]&\leq& C(\Lambda)^pm_p(e)^p
\EE \bigl[\EE \bigl(i_r^2(e) \rrvert r_S
\bigr)^p \bigr]
\\
&\leq&C(\Lambda)^{pL}m_p(e)^p \biggl(\sum
_{e'\in S}\EE \bigl[i_r^2
\bigl(e'\bigr) \bigr] \biggr)^p,
\end{eqnarray*}
where the last inequality follows from the second part of Lemma~\ref
{lem:regularite2}. Gathering this inequality and (\ref{eq:FSFe}),
\begin{eqnarray*}
\|\R_{S}\|_p&\leq& \bigl(2C(\Lambda)
\bigr)^{L}\min_{e\in S}m_p(e)\sum
_{e'\in S}\EE \bigl[i_r^2
\bigl(e'\bigr) \bigr]
\\
&\leq& \bigl(2C(\Lambda)\bigr)^{L}\sum_{e'\in S}m_p
\bigl(e'\bigr)\EE \bigl[i_r^2
\bigl(e'\bigr) \bigr]
\\
&\leq& C'(\Lambda)^L\sum
_{e'\in S}\bigl\|\R^{u,v}_{\{e'\}}\bigr\|_p,
\end{eqnarray*}
using the first part of the lemma.
\end{pf}

\section{Central limit theorem}
\label{sec:TLC}

Even if one takes two vertices $u$ and $v$ far apart, there is not
necessarily a Gaussian central limit theorem for the effective
resistance between them, since the influence of an edge near $u$, for
instance, may well represent a positive fraction of the total variance
of the resistance. However, let us define the \textit{influence of an
edge $e$} on the effective resistance between $u$ and $v$ by
\[
I^{u,v}(e):=\bigl\|\R^{u,v}_{\{e\}}\bigr\|_2^2.
\]
Then, if the maximal influence of an edge is small with respect to the
variance and if the graph has homogeneous currents and bounded degree,
one may obtain a Gaussian approximation. The following theorem shows a
result in this direction, whereas another instance of this phenomenon
will be described on a sequence of finite graphs, the discrete tori, in
Section~\ref{sec:tori}.

\begin{theo}
\label{theo:TLCinfini}
Let $\G$ be a graph with homogeneous currents and bounded degree,
equipped with elliptic resistances in $[1,\Lambda]$. For vertices $u$
and $v$ such that $\Var(\R^{u,v})>0$, define
\[
\beta(u,v)=\frac{\sup_{e\in\E}I^{u,v}(e)}{\Var(\R^{u,v})}
\]
and
\[
\overline{\R}^{u,v}:=\bigl[\R^{u,v}-\EE\bigl(\R^{u,v}
\bigr)\bigr]/{\sqrt{\Var\bigl(\R^{u,v}\bigr)}}.
\]
Let $\Phi$ be the standard Gaussian distribution function and let
$F^{u,v}$ be the distribution function of $\overline{\R}^{u,v}$. There
is a function $f\dvtx \RR^+\to[0,1]$, depending on $\G$ and $\Lambda$ only,
such that
\[
f(x)\mathop{\longrightarrow}_{x\to0}0
\]
and
\[
\bigl\|F^{u,v}-\Phi\bigr\|_\infty\leq f\bigl(\beta(u,v)\bigr).
\]
\end{theo}

\begin{pf}
Let us fix the vertices $u$ and $v$. For every integer $L$, one defines
$J=J(L)$ as
\[
J(L)=\bigl\{S\subset\E\mbox{ s.t. }\operatorname{diam}(S)\leq L\mbox{ and }S\neq\varnothing \bigr\}.
\]
Let $W_L$ be the random variable:
\[
W_L=\frac{\sum_{S\in J(L)}\R^{u,v}_S}{\sqrt{\sum_{S\in J(L)}\|\R
^{u,v}_S\|_2^2}}.
\]
Since $G$ has homogeneous currents, we know that for $L$ large\vspace*{1pt} enough,
$W_L$ will be close (in $L^2$-norm) to $\overline{\R}^{u,v}$. On the
other hand, since $\R^{u,v}_S$ depends only on $(r(e))_{e\in S}$, we
know that for every fixed $L$, $W_L$ is a sum of random variables with
only local dependence. Thus, one may use the work of \cite{ChenShao04}
to control the distance to normality.

To be more precise, let $F_L$ be the distribution function of $W_L$.
For $S$ in $J(L)$, define, using the notation of \cite{ChenShao04},
\begin{eqnarray*}
A_S &=& \bigl\{S_1\in J(L)\mbox{ s.t. }S\cap
S_1\neq\varnothing\bigr\},
\\
B_S &=& \bigl\{S_2\in J(L)\mbox{ s.t. }\exists
S_1\in A_S, S_1\cap S_2\neq
\varnothing\bigr\},
\\
C_S &=& \bigl\{S_3\in J(L)\mbox{ s.t. }\exists
S_2\in B_S, S_2\cap S_3\neq
\varnothing\bigr\}.
\end{eqnarray*}
Finally, let $N(C_S)=\{S'\in J(L)\mbox{ s.t. }C_S\cap C_{S'}\neq\varnothing\}$ and
\[
\kappa=\max_S\bigl\{\bigl|N(C_S)\bigr|,\bigl|\bigl
\{S'\mbox{ s.t. }S\in C_{S'}\bigr\}\bigr|\bigr\}.
\]
It is clear that
\[
C_S\subset\bigl\{S_3\in J(L)\mbox{ s.t. }\exists e\in
S,e'\in S_3, d\bigl(e,e'\bigr)\leq2L
\bigr\}.
\]
Thus,
\[
N(C_S)\subset\bigl\{S'\in J(L)\mbox{ s.t. }\exists e
\in S,e'\in S', d\bigl(e,e'\bigr)\leq6L\bigr\},
\]
which shows that
\[
\bigl|N(C_S)\bigr|\leq|S|\delta^{6L}2^{\delta^L}
\]
and
\[
\bigl|\bigl\{S'\mbox{ s.t. }S\in C_{S'}\bigr\}\bigr|\leq|S|
\delta^{2L}2^{\delta^L}.
\]
Thus, one may use Theorem~2.4 of \cite{ChenShao04} with $p=3$ and
$\kappa\leq2^{8\delta^L}$ where $\delta$ is a bound on the degrees in
$\G$. We obtain
\[
\|F_L-\Phi\|_{\infty}\leq\kappa\frac{\sum_{S\in J(L)}\llVert \R
^{u,v}_S\rrVert _3^3}{ (\sum_{S\in J(L)}\|\R^{u,v}_S\|_2^2
)^{3/2}}.
\]
Using Corollary~\ref{coro:firstlevel} and Lemma~\ref{lemm:regmean},
\begin{eqnarray*}
\|F_L-\Phi\|_{\infty}&\leq&C_1(L,\Lambda,\G)
\frac{\sum_{S\in J(L)}
(\sum_{e\in S}\|\R^{u,v}_{\{e\}}\|_3 )^3}{ (\sum_{e\in\E}\|\R
^{u,v}_{\{e\}}\|_2^2 )^{3/2}}
\\
&\leq&C_1(L,\Lambda,\G)\frac{\sum_{S\in J(L)}|S|^{2/3}\sum_{e\in S}\|\R
^{u,v}_{\{e\}}\|_3^3}{ (\sum_{e\in\E}\|\R^{u,v}_{\{e\}}\|_2^2
)^{3/2}}
\\
&\leq&C_1'(L,\Lambda,\G)\frac{\sum_{e\in\E}\|\R^{u,v}_{\{e\}}\|
_2^3}{ (\sum_{e\in\E}\|\R^{u,v}_{\{e\}}\|_2^2 )^{3/2}},
\end{eqnarray*}
one gets
\[
\|F_L-\Phi\|_{\infty}\leq C(L)\sqrt{\beta(u,v)},
\]
where $C(L)=C(\Lambda,\G,L)$ is a positive nondecreasing function of $L$.

On the other hand, Theorem~\ref{theo:smalldiam} ensures that
\[
\biggl\|\R^{u,v}-\EE\bigl(\R^{u,v}\bigr)-\sum
_{S\in J(L)}\R^{u,v}_S\biggr\|^2\leq
\eps(L)\Var \bigl(\R^{u,v}\bigr),
\]
where $\eps(L)=\eps(\Lambda,\G,L)$ is a positive nonincreasing
function of $L$ which goes to zero as $L$ goes to infinity. This implies
\[
\bigl\|\overline{\R}^{u,v}-W_L\bigr\|^2\leq4\eps(L).
\]
Notice that $\Phi$ is 1-Lipschitz (in fact, $1/\sqrt{2\pi}$-Lipschitz),
so for any $\eta>0$ and any $t\in\RR$
\begin{eqnarray*}
F^{u,v}(t)-\Phi(t)&\leq&F^{u,v}(t)-\Phi(t+\eta)+\eta
\\
&\leq&F^{u,v}(t)-F_L(t+\eta)+C(L)\sqrt{\beta(u,v)}+\eta
\\
&\leq&\PP\bigl(W_L-\overline{\R}^{u,v}>\eta\bigr)+C(L)
\sqrt{\beta(u,v)}+\eta
\\
&\leq&\frac{\|\overline{\R}^{u,v}-W_L\|_2^2}{\eta^2}+C(L)\sqrt{\beta (u,v)}+\eta
\\
&\leq&\frac{4\eps(L)}{\eta^2}+C(L)\sqrt{\beta(u,v)}+\eta.
\end{eqnarray*}
Symmetrically, one gets, for any $\eta>0$ and any $t\in\RR$
\[
\bigl|F^{u,v}(t)-\Phi(t)\bigr|\leq\frac{4\eps(L)}{\eta^2}+C(L)\sqrt{\beta (u,v)}+\eta.
\]
Optimizing in $\eta$ gives
\[
\bigl\|F^{u,v}-\Phi\bigr\|_{\infty}\leq3\eps^{1/3}(L)+C(L)\sqrt{
\beta(u,v)}.
\]
It remains to optimize in $L$. Let $L_0(x)=\sup\{L\in\NN\mbox{ s.t.
}3\eps^{1/3}(L)\geq\sqrt{x}C(L)\}$. Then
\[
\bigl\|F^{u,v}-\Phi\bigr\|_{\infty}\leq6\eps^{1/3}
\bigl(L_0\bigl(\beta(u,v)\bigr)\bigr).
\]
Since $L_0(x)$ goes to infinity as $x$ goes to zero, $f(x):=6\eps
^{1/3}(L_0(x))$ answers the theorem.
\end{pf}
It is in general difficult to apply this result because it is difficult
to bound $\beta(u,v)$. However, notice that the influence of an edge is
always bounded. Thus, on a bounded graph with homogeneous currents, if
the variance of $\R^{u,v}$ goes to infinity as $u$ and $v$ move apart,
one gets a central limit theorem. Notice that the last point is
equivalent to showing that $\EE[\sum_{e}(i_r^{u,v}(e))^4]$ goes to
infinity. It may be shown, for instance, that this is true on some wedges
of $\ZZ^2$, using the idea of Nash--Williams inequality. For instance,
let $h(x)=x^\alpha$ with $\alpha\leq1/3$, let $\V=\{(x,y)\in\ZZ^2\mbox{ s.t. }|y|\leq h(|x|)$ and let $\G$ be the subgraph of $\ZZ^2$ induced
by $\V$. Then one may derive a central limit theorem for $\R^{0,v}$ on
$\G$ when the distance $d(0,v)$ goes to infinity. Since this example is
not quasi-transitive, one has to use the Harnack inequality to prove
that the graph has homogeneous currents.

Notice also that already on $\ZZ^2$, the variance of the resistance is
only of order~1 (cf. \cite{NaddafSpencer}), and thus one cannot expect
a central limit theorem for point-to-point effective resistance when
the resistances are i.i.d. (since the influence of the edges near the
source and the sink is of order 1). In this respect, the interest of
Theorem~\ref{theo:TLCinfini} is rather limited. However, one should
rather think of it as a first step, with a clean statement, toward
central limit theorems for resistances on sequences of finite graphs,
as will be made clear in Section~\ref{sec:tori}.


\section{CLT for the effective resistance of the $d$-dimensional torus}
\label{sec:tori}

In this section, we investigate when $n$ becomes large the effective
resistance of the torus $\TT_n^d$ equipped with nonconstant i.i.d.
resistances from $[1,\Lambda]$. Here, $\TT_n^d$ is the graph $(\V_n^d,\E
_n^d)$ where $\V_n^d=(\ZZ/n\ZZ)^d$ and $\E_n^d$ is the set of oriented
edges of the torus: two vertices $x$ and $y$ of $\V_n^d$ are joined by
an edge from $x$ to $y$ if there is some $i\in\{1,\ldots,n\}$ such that
$x_i-y_i\in\{-1,1\}$ and for all $j\neq i$, $x_j=y_j$. One chooses also
exactly one edge of each orientation as follows:
\[
\E_{1/2}^n=\bigl\{(x,y)\in\E_n^d
\mbox{ s.t. }\exists i\in\{1,\ldots,n\} y_i-x_i=1\bigr\}.
\]
Recall that $\ell^2_-(\E_n^d,r)$ is the Hilbert space
\[
\ell^2_-\bigl(\E_n^d,r\bigr)= \bigl\{\theta
\in\RR^{\E_n^d}\mbox{ s.t. }\calE _r(\theta)<\infty\mbox{ and
}\forall e\in\E_n^d, \theta(e)=-\theta (-e) \bigr\},
\]
where
\[
\calE_r(\theta):=\sum_{e\in\E_{1/2}^n}r(e)
\theta^2(e),
\]
endowed with the scalar product
\[
\bigl(\theta,\theta'\bigr)_r=\sum
_{e\in\E_{1/2}^n}r(e)\theta(e)\theta'(e).
\]
Also, for resistances $r$ in $[1,\Lambda]^{\E_{1/2}^n}$, all the sets
$\ell^2_-(\E_n^d,r)$ are the same, and we denote this space by $\ell
^2_-(\E_n^d)$.

\begin{figure}

\includegraphics{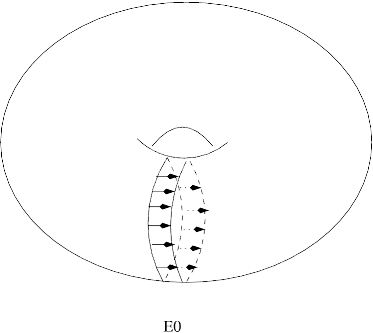}

\caption{The cut $E_0$ on a 2-dimensional torus.}
\label{fig:torusE0}
\end{figure}

Since $\TT_n^d$ has no boundary, our first objective is to define the
effective resistance in a natural and translation invariant way. First,
we define a special \textit{cut along direction} 1 (see Figure~\ref{fig:torusE0}):
\[
E_0:=\bigl\{(x,y)\in\E_{1/2}^n \mbox{ s.t. }x
\sim y, x_1=0\mbox{ and }y_1=1\bigr\},
\]
and the flows which \textit{cross the torus along direction} 1, \textit{with
intensity} 1:
\[
\Theta_0^n:=\biggl\{\theta\in\ell^2_-
\bigl(\E_n^d\bigr)\mbox{ s.t. }d^*\theta=0\mbox{ and }\sum
_{e\in E_0}\theta(e)=1\biggr\}.
\]

Notice that the elements of $\Theta_0^n$ are sourceless flows and that
the definition is independent of the choice of the cut (along direction
1). Indeed, if we define, for any $i$ in $\{0,\ldots,n-1\}$,
\[
E_i:=\bigl\{(x,y)\in\V_n^d\times
\V_n^d\mbox{ s.t. }x\sim y, x_1=i\mbox{ and
}y_1=i+1\bigr\},
\]
then, for any $i$ and any $\theta\in\ell^2_-(\E_{1/2}^n)$, we have
%
\begin{equation}
\label{eq:EiEiplus1}
\sum_{e\in E_{i}}\theta(e)-\sum
_{e\in E_{i-1}}\theta(e)=\sum_{x
\, \,\mathrm{s.t.} \,\,  x_1=i}d^*
\theta(x).
\end{equation}
Thus, for any sourceless flow $\theta$, $\sum_{e\in E_0}\theta(e)$
equals 1 if and only if $\sum_{e\in E_i}\theta(e)=1$ for some $i$ in $\{
0,\ldots,n-1\}$. Thus, for any translation $\tau$ on the torus,
%
\begin{equation}
\label{eq:translationTheta0}
\theta\in\Theta_0^n \quad \Longleftrightarrow \quad
\theta\circ\tau\in\Theta_0^n.
\end{equation}

\begin{defi}
One defines \textit{the effective resistance of the torus} as
%
\begin{equation}
\label{eq:resistancetorus}
\R_n:=\inf_{\theta\in\Theta_0^n}
\calE_r(\theta).
\end{equation}
\end{defi}

Using the Nash--Williams inequality, it is easy to show that $\R
_n=\Theta(1/n^{d-2})$. Using the ideas of the preceding sections, one
may show that $\R_n$ satisfies a central limit theorem. This is the
main result of the present paper.

\begin{theo}
\label{theo:CLTTorus}
Suppose that $r(e)$, $e\in\E_{1/2}^n$ are i.i.d. with positive variance
and support in $[1,\Lambda]$ for some $\Lambda>1$. Then
\[
\Var(\R_n)=\Theta\bigl(n^{4-3d}\bigr),
\]
and the resistance satisfies a central limit theorem:
\[
\frac{\R_n-\EE(\R_n)}{\sqrt{\Var(\R_n)}}\mathop{\longrightarrow}_{n\rightarrow\infty}^{\L}\N(0,1).
\]
\end{theo}

It is straightforward to show that a similar statement holds for what
should be called \textit{the effective conductance of the torus $\TT
_n^d$}, $\C_n:=1/\R_n$ or \textit{the mean conductivity of the torus $\TT
_n^d$}: $A_n=n^{2-d}\C_n$. One obtains
\[
\Var(A_n)=\Theta\bigl(n^{-d}\bigr)
\]
and
\[
\frac{A_n-\EE(A_n)}{\sqrt{\Var(A_n)}}\mathop{\longrightarrow}_{n\to\infty}^{\L}\N (0,1).
\]
In \cite{Boivin09}, the definition of $A_n$ is different, based on the
discrete cube and not on the torus, however the behaviour should be the
same. \cite{Boivin09} obtains only suboptimal bounds for the variance
of the mean conductivity. In \cite{GloriaOttoVariance}, an optimal
bound on the variance is obtained for a related quantity based on the
\textit{corrector} of homogenization theory. Again, the behaviour should
be the same. In any case, the central limit theorem is new.

Our strategy to show Theorem~\ref{theo:CLTTorus} is simply to apply the
ideas of Theorem~\ref{theo:TLCinfini} to this setting, where the
infinite graph is traded against a growing sequence of finite graphs.
First, notice that there is a unique minimal flow reaching the infimum
in the definition of the resistance. It is characterized by an
orthogonality criterion which is the analog of the Kirchhoff cycle law
(see Definition~\ref{defi:Kirchhoff}) in this setting, so we shall call
it a \textit{pseudo-Kirchhoff cycle law}. Define the following tangent
vector space to $\Theta_0^n$:
\[
\overrightarrow{\Theta}:=\biggl\{\theta\in\ell^2_-\bigl(
\E_n^d\bigr)\mbox{ s.t. }d^*\theta=0\mbox{ and }\sum
_{e\in E_0}\theta(e)=0\biggr\}.
\]

\begin{lemm}
\label{lemm:pseudoKirchhoff}
The infimum in (\ref{eq:resistancetorus}) is attained at a unique flow
$i_{r,n}^{\mathrm{per}}$ which is the unique flow in $\Theta_0^n$ that satisfies
the following pseudo-Kirchhoff cycle law:
\[
\forall h\in\overrightarrow{\Theta},\qquad \bigl(h,i_{r,n}^{\mathrm{per}}
\bigr)_r=0.
\]
\end{lemm}

\begin{pf}
This is standard since $(\cdot,\cdot)_r$ is a scalar product. Let $\theta_0\in
\Theta_0^n$ be fixed. Notice that $\Theta_0^n=\theta_0+ \overrightarrow
{\Theta}$. Define $i_{r,n}^{\mathrm{per}}$ as the orthogonal projection of
$\theta_0$ on $\overrightarrow{\Theta}^{\perp}$ for $(\cdot,\cdot)_r$. It is
the unique element $ i $ of $\ell_-^2(\E_n^d)$ satisfying
\[
\theta_0- i \in\overrightarrow{\Theta} \quad\mbox{and}\quad \forall h\in
\overrightarrow{\Theta},\qquad (h, i )_r=0,
\]
which is equivalent to
\[
i \in\Theta_0^n\quad \mbox{and}\quad \forall h\in
\overrightarrow{\Theta }, \qquad (h, i )_r=0.
\]
This shows the last part of the lemma. Now, for any $\theta\in\Theta
_0^n$, $\theta-i_{r,n}^{\mathrm{per}}$ belongs to~$\overrightarrow{\Theta}$, and thus
\begin{eqnarray*}
\calE_r(\theta)&=&\calE_r\bigl(\theta-i_{r,n}^{\mathrm{per}}
\bigr)+2\bigl(\theta -i_{r,n}^{\mathrm{per}},i_{r,n}^{\mathrm{per}}
\bigr)_r+\calE_r\bigl(i_{r,n}^{\mathrm{per}}
\bigr)
\\
&=&\calE_r\bigl(\theta-i_{r,n}^{\mathrm{per}}\bigr)+
\calE_r\bigl(i_{r,n}^{\mathrm{per}}\bigr),
\end{eqnarray*}
which shows that the infimum of $\calE_r(\theta)$ over $\theta\in\Theta
_0^n$ is uniquely attained at $i_{r,n}^{\mathrm{per}}$.
\end{pf}
Since the minimal flow is sourceless but satisfies the additional
condition that the net flow through $E_0$ is 1, one needs to adapt the
setting of the first sections in order to prove
Theorem~\ref{theo:CLTTorus}. For any $e\in\E_{1/2}^n$, let
\[
\Theta_e^n:=\biggl\{\mbox{unit flows }\theta\mbox{ from
}e_-\mbox{ to }e_+\mbox{ s.t. }\sum_{e'\in E_0}\theta
\bigl(e'\bigr)=\II_{e\in E_0}\biggr\}.
\]
From (\ref{eq:EiEiplus1}), one sees that for any unit flow $\theta$
from $e_-$ to $e_+$, and any $i$,
\[
\sum_{e'\in E_i}\theta\bigl(e'\bigr)=\sum
_{e'\in E_0}\theta\bigl(e'\bigr)+
\II_{e\in E_i}-\II _{e\in E_0}.
\]
Thus, one sees that for any translation $\tau$ on the torus, and any
edge $e$,
%
\begin{equation}
\label{eq:translationThetae}
\theta\in\Theta_e^n \quad \Longleftrightarrow \quad\theta
\circ\tau\in\Theta_{\tau
^{-1}(e)}^n.
\end{equation}
Then, one may show as in the proof of Lemma~\ref{lemm:pseudoKirchhoff}
that $\theta\mapsto\calE_r(\theta)$ has a unique minimizer on $\Theta
_e^n$, that we shall call $j_{r,n}^e$, and which is characterized by
the same pseudo-Kirchhoff cycle law.

\begin{lemm}
\label{lemm:pseudoKirchhoffjrne}
The following infimum
\[
\inf_{\theta\in\Theta_e^n}\calE_r(\theta)
\]
is attained at a unique flow $j_{r,n}^e$, which is the orthogonal
projection of $\chi_e$ on $\overrightarrow{\Theta}^{\perp}$ in $\ell
_-^2(\E_n^d,r)$. It is the unique flow in $\Theta_e^n$ that satisfies
the pseudo-Kirchhoff cycle law:
\[
\forall h\in\overrightarrow{\Theta}, \qquad\bigl(h,j_{r,n}^e
\bigr)_r=0.
\]
\end{lemm}

\begin{pf}
The proof is completely similar to the proof of Lemma~\ref
{lemm:pseudoKirchhoff}, since $\overrightarrow{\Theta}$ is again the
tangent vector space to $\Theta_e^n$.
\end{pf}
The role of $j_{r,n}^e$ will be similar to that of $i_r^e$ in the first
sections, as hinted by the following lemma.

\begin{lemm}
\label{lem:regularitytorus}
The functions $r\mapsto i_{r,n}^{\mathrm{per}}(e)$, for any edge $e$, and
$r\mapsto\R_n(r)$ admit partial derivatives of all orders. In
addition, for any edges $e$, $e'$:
\begin{longlist}[(iii)]
\item[(i)] $\forall e'\neq e, \partial_{e'}
i_{r,n}^{\mathrm{per}}(e)=\frac{i_{r,n}^{\mathrm{per}}(e')}{r(e')}j^{e'}_{r,n}(e)=\frac
{i_{r,n}^{\mathrm{per}}(e')}{r(e)}j^{e}_{r,n}(e') $.

\item[(ii)] $\forall e, \partial_e
i_{r,n}^{\mathrm{per}}(e)=\frac{i_{r,n}^{\mathrm{per}}(e)}{r(e)}(j^{e}_{r,n}(e)-1) $.\vspace*{2pt}
\item[(iii)] $\forall e, \partial_e \R
_n(r)=(i_{r,n}^{\mathrm{per}}(e))^2 $.
\end{longlist}
\end{lemm}

\begin{pf}
The fact that $r\mapsto i_{r,n}^{\mathrm{per}}(e)$ and $r\mapsto\R_n(r)$ admits
partial derivatives of all order is analogous to the classical case;
cf. the proof of Lemma~\ref{lem:regularite1}.

Let us fix some edge $e'\in\E_{1/2}^n$. One may thus differentiate the
node law and the pseudo-Kirchhoff cycle law of Lemma~\ref
{lemm:pseudoKirchhoff} with respect to $r(e')$ to obtain
\[
\forall x\in\V_n^d, \qquad  d^*\bigl[\partial_{e'}
i_{r,n}^{\mathrm{per}}\bigr](x)=0
\]
and
\[
\forall h\in\overrightarrow{\Theta},\qquad \sum_{e}h(e)r(e)
\biggl(\partial_{e'} i_{r,n}^{\mathrm{per}}(e)+
\frac{i_{r,n}^{\mathrm{per}}(e')}{r(e')}\chi _{e'}(e) \biggr)=0.
\]
Thus, if we define
\[
\theta(e)=\partial_{e'} i_{r,n}^{\mathrm{per}}(e)+
\frac
{i_{r,n}^{\mathrm{per}}(e')}{r(e')}\chi_{e'}(e),
\]
we see that
\[
\forall x\notin e',\qquad  d^*\theta(x)=0
\]
and
\[
\forall h\in\overrightarrow{\Theta},\qquad (h,\theta)_r=0.
\]
Furthermore,
\[
d^*\theta\bigl(e'_-\bigr)=\frac{i_{r,n}^{\mathrm{per}}(e')}{r(e')}
\]
and
\[
\sum_{e\in E_0}\theta(e)=\frac{i_{r,n}^{\mathrm{per}}(e')}{r(e')}
\II_{e'\in
E_0}.
\]
It follows from the characterization of $j_{r,n}^{e'}$ that
\[
\theta=\frac{i_{r,n}^{\mathrm{per}}(e')}{r(e')}j_{r,n}^{e'}.
\]
This gives the proof of the first two equations. The proof of the last
one is analogous to the classical case; cf. Lemma~\ref{lem:regularite1}.
\end{pf}
When $n$ is large, we would like to compare $j_{r,n}^e$ to a flow on
the whole lattice $\ZZ^d$. To do this, we shall couple the network $(\ZZ
^d,\EE^d)$ with all the tori by ``unwrapping'' each torus on $\ZZ^d$.
\textit{This construction will be used throughout the section.} Since our
main objects are elements of $\Theta_0^n$ and $\Theta_e^n$, we mainly
need to identify the set of edges $\E_n^d$, equipped with their
resistances, as subsets of $\EE^d$, and then be careful about what
happens to the boundary operator through this identification.

First, fix $e$ to be any edge such that $e_-$ is the origin of $\ZZ^d$.
Let $r\in[1,\Lambda]^{\EE_{1/2}^d}$ be a fixed collection of
resistances and define
\[
\V_n=\bigl\{-\lfloor n/2\rfloor,\ldots,\bigl\lfloor(n-1)/2\bigr
\rfloor\bigr\}^d,
\]
where $\lfloor\cdot\rfloor$ is the integer part, and
\[
\E_{1/2}^n=\bigl\{(x,y)\in\V_n\times
\ZZ^d\mbox{ s.t. }\exists i\in\{1,\ldots,n\} y_i-x_i=1
\mbox{ and }\forall j\neq i, x_j=y_j\bigr\}
\]
so that these sets are roughly centered around the origin. Now, we let
$\G_n^{\mathrm{per}}$ be the network with edge set $\E_{1/2}^n$, resistances
induced by $r$ and set of vertices all the endpoints of the edges of $\E
_{1/2}^n$ with periodic condition, that is, vertices with identical
coordinates modulo $n$ are identified (this takes care of the boundary
operator on the torus). Clearly, $\G_n^{\mathrm{per}}$ is isomorphic to $\TT
_n^d$, and we shall thus use the notation $E_0$, $\Theta_e^n$, $\Theta
_0^n$ and $j_{r,n}^e$ on $\G_n^{\mathrm{per}}$ as well.

\begin{lemm}
\label{lemm:convergencejrne}
Let $e$ be any fixed edge such that $e_-$ is the origin of $\ZZ^d$. Let
$i_r^e$ be the minimal current on $\ZZ^d$ from $e_-$ to $e_+$, then
$j_{r,n}^e$ converges to $i_r^e$ in $\ell^2_-(\EE^d,r)$.
\end{lemm}

\begin{pf}
We let $\G_n^F$ be the subgraph of $\ZZ^d$ induced by the set of
vertices $\V_n$. Also, we define $\G_n^W$ be the graph obtained from $\ZZ
^d$ by identifying all the vertices outside $\{-\lfloor n/2\rfloor
+1,\ldots,\lfloor(n-1)/2\rfloor\}^d$. See Figure~\ref{fig:GnGnFGnW}.
We equip these graphs with resistances given by $r$. Since these graphs
have sets of edges which are still subsets of $\EE^d$, there is no
ambiguity about what resistance is assigned to which edge.
\begin{figure}

\includegraphics{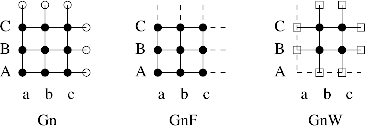}

\caption{The three graphs $\G_n^{\mathrm{per}}$, $\G_n^F$ and $\G_n^W$ when $n=3$
and $d=2$. Squares are identified as a single vertex, while empty
circles are subject to periodic identification modulo $n$.}
\label{fig:GnGnFGnW}
\end{figure}

We need now to introduce some terminology from the theory of electrical
networks. Let $\G=(\V,\E)$ be a graph equipped with resistances
$r=(r(e))_{e\in\E}$ and
suppose that $(\H_n)_{n\geq0}$ is a sequence of finite subgraphs of
$\G$ that exhausts $\G$, as in Section~\ref{sec:prel_current}. A
\textit{star} of a graph $\G=(\V,\E)$ is a member of $\ell^2_-(\E,r)$ of the
form $\sum_{e_-=x}c(e)\chi_e$ for $x\in\V$. Let $\bigstar$ (resp.,
$\bigstar_n$) denote the closed subspace spanned by~the stars in $\ell
^2_-(\E,r)$ (resp., by the stars of $\H_n^W$). A \textit{cycle} is a
member of $\ell^2_-(\E,r)$ of the form $\sum_{i=1}^n\chi^{e_k}$ with
$e_1,\ldots,e_n$ an oriented cycle in $\G$. Let $\diamondsuit$ (resp.,
$\diamondsuit_n$) denote the closed subspace spanned by the cycles in
$\ell^2_-(\E,r)$ (resp., by the cycles of $\H_n$). To understand the
introduction of wired and free networks, notice that all stars of the
wired network $\H_n^W$ are stars in $\G$, except for the star at the
``extra vertex'' which represents all the vertices outside $\H_n$, but
this additional star is just the opposite of all the stars in $\H_n$,
thus $\bigstar_n$ is a subspace of $\bigstar$. Furthermore, since $(\H
_n)_n$ exhausts $\G$, each star of $\G$ is a member of $\bigstar_n$ for
$n$ large enough. This shows that $\bigstar=\overline{\bigcup_n \bigstar
_n}$. Also, all the cycles of $\H_n$ are cycles of $\G$, and thus
$\diamondsuit_n$ is a subspace of $\diamondsuit$, and each cycle of
$\diamondsuit$ is in $\diamondsuit_n$ for $n$ large enough. This shows
that $\diamondsuit=\overline{\bigcup_n \diamondsuit_n}$.

For a closed subspace $V$, denote by $P_V$ the projection on $V$ in
$\ell^2_-(\E,r)$. Then, for any edge $e$ in $\H_0$, $P_{\bigstar_n}\chi
_e$ is the unique current on $\H_n^W$ between $e_-$ and $e_+$ and it
converges in $\ell^2_-(\E,r)$ to $i_r^{e}=P_\bigstar\chi_e$, the
minimal current on $\G$ from $e_-$ to $e_+$ (this is a simple
consequence of the fact that $\bigstar=\overline{\bigcup_n \bigstar_n}$;
see Exercise 9.2 in~\cite{LyonsPeres2011}, and also Propositions 9.1
and 9.2 therein). Also, $P_{\diamondsuit_n^{\perp}}\chi_e$ converges to
$i_r^{F,e}=P_{\diamondsuit^{\perp}}\chi_e$, the \textit{free current}
from $e_-$ to $e_+$. In general, the free current $i_r^{F,e}$ is a
current that may or may not be equal to $i_r^{e}$. However, on $\ZZ^d$
(and on finite graphs of course), it is known that those currents
coincide. One says in this case that ``currents are unique.'' Another
way to express this is to say that $\bigstar=\diamondsuit^{\perp}$. An
argument can be given as follows: this is trivial on $\ZZ$, because
free and wired currents between $u$ and $v$ can be easily computed to
be just 1 from $u$ to $v$ and 0 elsewhere, then the unicity of currents
is preserved under cartesian product (Exercise 9.7 in~\cite
{LyonsPeres2011}), this shows currents are unique on $\ZZ^d$ with unit
resistances. Finally, the unicity of currents is preserved under
``rough isometries'' (Theorem~9.9 in~\cite{LyonsPeres2011}), which
include elliptic perturbations of the weights.

Now, we return to our setting and let $\bigstar_n$ (resp., $\bigstar
_n^{\mathrm{per}}$) be the linear subspace of $\ell^2_-(\E_{1/2}^n)$ spanned by
the stars of $\G_n^W$ (resp., $\G_n^{\mathrm{per}}$). Let also $\diamondsuit_n$
(resp., $\diamondsuit_n^{\mathrm{per}}$) the linear subspace\vspace*{1pt} spanned by the
cycles of $\G_n^F$ (resp., $\G_n^{\mathrm{per}}$). Since $\G_n^F$ is a strict
subgraph of $\G_n^{\mathrm{per}}$,
\[
\diamondsuit_n\subset\diamondsuit_n^{\mathrm{per}}.
\]
Furthermore, any cycle in $\G_n^F$ must traverse $E_0$ in one direction
the same number of times that it traverses it in the other direction
(notice that this is not true on $\G_n^{\mathrm{per}}$, due to the periodic
boundary conditions). Thus,
\[
\diamondsuit_n\subset\diamondsuit_n^{\mathrm{per}}\cap
\biggl\{\theta\in\ell_-^2\bigl(\E _{1/2}^n\bigr)
\mbox{ s.t. }\sum_{e\in E_0}\theta(e)=0\biggr\}.
\]
Also, note that the stars in $\G_n^W$ are generated by the stars at
vertices in $\{-\lfloor n/2\rfloor+1,\ldots,\lfloor(n-1)/2\rfloor\}
^d$, since the star at the ``exterior vertex'' equals the opposite of
the sum of all the other stars. But all those stars are stars of
$\G_n^{\mathrm{per}}$. Thus,
\[
\bigstar_n\subset\bigstar_n^{\mathrm{per}}.
\]
Recall from Lemma~\ref{lemm:pseudoKirchhoffjrne} that $j_{r,n}^e$ is
the orthogonal projection in $\ell_-^2(\E_n^d,r)$ of $\chi_e$ on
$H_n=\overrightarrow{\Theta}^{\perp}$. It is easy to see that
\[
\overrightarrow{\Theta}=\biggl\{\theta\in\diamondsuit_n^{\mathrm{per}}
\mbox{ s.t. }\sum_{e\in E_0}\theta(e)=0\biggr\}.
\]
Indeed, on $\G_n^{\mathrm{per}}$, the condition $d^*\theta=0$ may be written as
$\theta\in(\bigstar_n^{\mathrm{per}})^{\perp}$, which is equivalent to $\theta
\in\diamondsuit_n^{\mathrm{per}}$ since $\G_n^{\mathrm{per}}$ is a finite graph. Notice
that $j_{r,n}^e$ is a flow between $e_-$ and $e_+$ on $\G_n^{\mathrm{per}}$ but
not necessarily on $\ZZ^d$. The inclusions above give
\[
\bigstar_n\subset H_n \subset\diamondsuit_n^{\perp}.
\]
Let $\bigstar$ (resp., $\diamondsuit$) denote the closed linear
subspace spanned by stars (resp., cycles) in $\ell^2_-(\EE^d,r)$. Then,
according to the elements of the theory of electrical networks stated
above, $P_{\bigstar_n}\chi_e$ (resp., $P_{\diamondsuit_n^{\perp}}\chi
_e$) converges to $i_r^{e}=P_\bigstar\chi_e$, the minimal current on
$\ZZ^d$ from $e_-$ to $e_+$ (resp., $i_r^{F,e}=P_{\diamondsuit^{\perp
}}\chi_e$, the free current on $\ZZ^d$ from $e_-$ to $e_+$). But on $\ZZ
^d$, currents are unique. As a consequence, $(j_{r,n}^e)_{n\geq1}$
converges also to $i_r^{e}$, the minimal current on $\ZZ^d$ from $e_-$
to $e_+$.
\end{pf}

In order to adapt the notion of graphs with homogeneous currents to
this setting of sequences of finite graphs, define
\[
\alpha(d,L,\Lambda):=\sup_{n\geq1}\sup_{r\in[1,\Lambda]^{\E_{1/2}}}
\sup_{e\in\TT_n^d}\sum_{d(e',e)\geq L}r
\bigl(e'\bigr) \bigl(j_{r,n}^e
\bigl(e'\bigr)\bigr)^2.
\]

\begin{prop}
\label{prop:homogeneoustorus}
The $d$-dimensional discrete tori have homogeneous currents in the
sense that for any $\Lambda\geq1$,
\[
\alpha(d,L,\Lambda)\mathop{\longrightarrow}_{L\rightarrow+\infty} 0.
\]
\end{prop}

\begin{pf}
We will adapt the proof of Proposition~\ref{prop:Mikael}, using the
convergence of~$j_{r,n}^e$ given by Lemma~\ref{lemm:convergencejrne}.
Let $\V=\ZZ^d$ and $\E=\EE^d$.

Let $e_1,\ldots, e_d$ be the $d$ edges going from $0$ to some point of
nonnegative coordinates. Thanks to the translation property (\ref
{eq:translationThetae}), one has, for any translation $\tau$ on the torus:
\[
j^{\tau(e)}_{r\circ\tau^{-1},n}\circ\tau=j^{e}_{r,n}.
\]
Thus,
\[
\alpha(d,L,\Lambda):=\sup_{n\geq1}\sup_{e\in\{e_1,\ldots,e_d\}}
\sup_{r\in[1,\Lambda]^{\E_{1/2}}}\sum_{d(e',e)\geq
L}r
\bigl(e'\bigr) \bigl(j_{r,n}^e
\bigl(e'\bigr)\bigr)^2.
\]
Let $e$ be any fixed edge such that $e_-$ is the origin, and define
\[
c:=\limsup_{L\rightarrow\infty}\sup_{n\geq1}\sup
_{r\in[1,\Lambda]^{\E
_{1/2}}}\sum_{d(e',e)\geq L}r
\bigl(e'\bigr) \bigl(j_{r,n}^e
\bigl(e'\bigr)\bigr)^2.
\]
It is thus enough to prove that $c=0$. Notice that when $n$ is fixed,
\[
\limsup_{L\rightarrow\infty}\sup_{r\in[1,\Lambda]^{\E_{1/2}}}\sum
_{d(e',e)\geq L}r\bigl(e'\bigr) \bigl(j_{r,n}^e
\bigl(e'\bigr)\bigr)^2=0,
\]
since the sequence in $L$ is zero for $L$ large enough. Thus, one may
find a sequence $(r_L,n_L)_{L\geq1}$ in $[1,\Lambda]^{\E_{1/2}}\times
\NN$ such that
\[
c=\lim_{L\to\infty}\sum_{d(e',e)\geq L}r_L
\bigl(e'\bigr) \bigl(j_{r_L,n_L}^e
\bigl(e'\bigr)\bigr)^2,
\]
and $n_L\ds\mathop{\longrightarrow}_{L\to\infty}+\infty$. The sequence
$j_{r_L,n_L}^e$ is bounded in $\ell^2_-(\EE^d)$. Thus, by compactness
of $[1,\Lambda]^{\E_{1/2}}$ one may extract a sequence from
$(r_L,n_L)_{L\geq1}$, that we shall still denote by $(r_L,n_L)_{L\geq
1}$ to lighten the notation, such\vspace*{1pt} that $(j_{r_L,n_L}^e(e'))_L$
converges $\theta(e')$ for any $e'$, and $(r_L(e'))_L$ converges to
some $r(e')\in[1,\Lambda]$ for any $e'$. Notice that $\theta$ is a
unit flow on the whole lattice $\ZZ^d$ since for any $L$,
$(j_{r_L,n_L}^e)_L$ is a unit flow on $G_n^{\mathrm{per}}$. Using the minimality
of $i_r^e$,
%
\begin{eqnarray}
\nonumber
\calE_r\bigl(i_r^e\bigr)&\leq&
\calE_r(\theta)
\\
\nonumber
&=&\sum_{e'\in\E_{1/2}}r\bigl(e'\bigr)
\theta^2\bigl(e'\bigr)
\\
\nonumber
&\leq&\limsup_{L\to\infty} \biggl(\sum
_{d(e',e)<L}r_L\bigl(e'\bigr)
\bigl(j_{r_L,n_L}^e\bigl(e'\bigr)
\bigr)^2 \biggr)
\\[-8pt]
\label{eq:majoMikaeltorus}
\\[-8pt]
\nonumber
&=&\limsup_{L\to\infty} \biggl(\calE_{r_L}
\bigl(j_{r_L,n_L}^e\bigr)-\sum_{d(e',e)\geq L}r_L
\bigl(e'\bigr) \bigl(j_{r_L,n_L}^e
\bigl(e'\bigr) \bigr)^2 \biggr)
\\
\nonumber
&= &\limsup_{L\to\infty}\calE_{r_L}
\bigl(j_{r_L,n_L}^e\bigr)-c
\\
\nonumber
 &\leq&\limsup_{L\to\infty}\calE _{r_L}
\bigl(j_{r,n_L}^e\bigr)-c,
\end{eqnarray}
where in the last inequality we used the minimality property of
$j_{r_L,n_L}^e$. Now, using Minkowski's inequality,
\begin{eqnarray*}
\calE_{r_L}\bigl(j_{r,n_L}^e\bigr)&=&\sum
_{e'\in\E
_{1/2}}r_L\bigl(e'\bigr)
\bigl(j_{r,n_L}^e\bigl(e'\bigr)
\bigr)^2
\\
&\leq& \biggl(\sqrt{\sum_{e'\in\E_{1/2}}r_L
\bigl(e'\bigr) \bigl(i_{r}^e
\bigl(e'\bigr)\bigr)^2}+\sqrt {\Lambda\sum
_{e'\in\E_{1/2}}\bigl(j_{r,n_L}^e
\bigl(e'\bigr)-i_{r}^e\bigl(e'
\bigr)\bigr)^2} \biggr)^2.
\end{eqnarray*}
Since $(r_L)_{L\geq1}$ converges simply to $r$ and is bounded by
$\Lambda$, the dominated convergence theorem gives
\[
\sum_{e'\in\E_{1/2}}r_L\bigl(e'
\bigr) \bigl(i_{r}^e\bigl(e'\bigr)
\bigr)^2 \mathop{\longrightarrow}_{L\to\infty}\calE_r
\bigl(i_r^e\bigr),
\]
Lemma~\ref{lemm:convergencejrne} says that $(j_{r,n_L}^e)_{L\geq1}$
converges to $i_r^e$ in $\ell^2_-(\E_{1/2})$. Thus,
\[
\limsup_{L\to\infty}\calE_{r_L}\bigl(j_{r,n_L}^e
\bigr)\leq\calE_r\bigl(i_r^e\bigr).
\]
Plugging this into (\ref{eq:majoMikaeltorus}) shows that $c=0$.
\end{pf}
%
Now, one may complete the proof of Theorem~\ref{theo:CLTTorus}.
\begin{pf*}{Proof of Theorem~\protect\ref{theo:CLTTorus}}
With Lemma~\ref{lem:regularitytorus} at hand, it is easy to reproduce
the proofs of Lemma~\ref{lem:regularite2}, Corollary~\ref
{coro:firstlevel} and Theorem~\ref{theo:smalldiam} with $\alpha(\G
,L,\Lambda)$ replaced by $\alpha(d,L,\Lambda)$. One obtains notably the
existence of constants $C(\Lambda)$ and $C(\Lambda,d)$ such that for
any $n$, and any $L\geq1$,
\[
\sum_{\operatorname{diam}(S)\geq L}\bigl\|(\R_n)_S
\bigr\|_2^2\leq C(\Lambda)\alpha(G,L,\Lambda )\sum
_{S\neq\varnothing}\bigl\|(\R_n)_S
\bigr\|_2^2
\]
and
\[
\sum_{e}\bigl\|(\R_n)_{\{e\}}
\bigr\|_2^2\leq\Var(\R_n)=\sum
_{S\neq\varnothing}\bigl\| (\R_n)_S
\bigr\|_2^2\leq C(\Lambda,d)\sum_{e}
\bigl\|(\R_n)_{\{e\}}\bigr\|_2^2.
\]
Now, thanks to the translation invariance of the model given by (\ref
{eq:translationTheta0}) and the fact that the edge-resistances are i.i.d.,
\[
\beta_n:=\sup_n\sup_{e\in\E_{1/2}^n}
\frac{\|(\R_n)_{\{e\}}\|_2^2}{\Var
(\R_n)}=\Theta\bigl(1/n^d\bigr)
\]
and
\[
\sup_{e\in\E_{1/2}^n}\EE\bigl[\bigl(i_{r,n}^{\mathrm{per}}(e)
\bigr)^2\bigr]=\Theta\biggl(\frac{1}{n^d}\EE (\R_n)
\biggr)=\Theta\biggl(\frac{1}{n^{2d-2}}\biggr).
\]
This already shows that
\[
\Var(\R_n)=\Theta\Bigl(\sup_{e\in\E_{1/2}^n}\EE \bigl[
\bigl(i_{r,n}^{\mathrm{per}}(e)\bigr)^2\bigr]^2
\Bigr)=\Theta\biggl(\frac{1}{n^{3d-4}}\biggr).
\]
Now, Proposition~\ref{prop:homogeneoustorus} allows to repeat the proof
of Theorem~\ref{theo:TLCinfini}. Let $F_n$ be the distribution function
of $\frac{\R_n-\EE(\R_n)}{\sqrt{\Var(\R_n)}}$. One obtains the
existence of a function $f$ having limit 0 at $0^+$ and such that for
any $n$,
\[
\|F_n-\Phi\|_{\infty}\leq f(\beta_n).
\]
This completes the proof of the central limit theorem.
\end{pf*}
%

\section{Perspectives}

We end this article with some questions left open.

First, it is not clear whether the notion of homogeneous currents is
really useful to get a central limit theorem. For instance, in the
counterexample of Figure~\ref{fig:nonhomogeneous}, one sees that the
currents $i_r^e$ still spread most of their mass at very localized
places, namely near $e$ and near the edge $e'_k$ which is in the same
connected component as~$e$. Thus, one could be able to adapt the proof
of the central limit theorem in this special case. One may wonder
whether the sole hypotheses of bounded degree and small influences are
enough to get a central limit theorem. On the other hand, if the
homogeneous currents hypothesis is proved really necessary, it would be
important to understand which graphs satisfy it, and whether it is
stable under perturbations like quasi-isometries.

Second, the most obvious question left open is the one raised in the
\hyperref[sec1]{Introduction}, that is to determine the order of the variance and to
show a central limit theorem for the resistance on the cube of side
length $n$ in $\ZZ^d$, and not only on the torus. More generally,
consider a domain $\Omega$ of $\RR^d$ with two disjoint\vspace*{1pt} subsets of its
boundary, $A$ and $Z$. Let $\G_n$ be the graph induced by $\Omega\cap
\frac{1}{n}\ZZ^d$ and let $\R_n$ be the effective resistance between
$A$ and $Z$ on $\G_n$. Then we conjecture that a Gaussian central limit
theorem holds for $\R_n$.

\section*{Acknowledgements}
I would like to warmly thank Itai Benjamini and Michel Benaim with whom
I began this research. Thanks also to Pierre Mathieu, Alano Ancona,
Daniel Boivin and Thierry Delmotte for very useful discussions and to
an anonymous referee for a very careful reading which helped to
improve this article. Last but not least, the main part of Proposition~\ref{prop:Mikael} is due to Mika\"el de la Salle; I owe a lot to him.






\printaddresses
\end{document}